\documentclass{article}
\usepackage{amsmath,amsthm,mathrsfs,amsfonts,amssymb,stmaryrd,xy,bm}
\usepackage[a4paper]{geometry}
\theoremstyle{definition}
\newtheorem{remark}{Remark}[section]
\newtheorem{remarks}[remark]{Remarks}

 \theoremstyle{plain}
\newtheorem{definition}[remark]{Definition}
\newtheorem{theorem}[remark]{Theorem}
\newtheorem{proposition}[remark]{Proposition}
\newtheorem{corollary}[remark]{Corollary}
\newtheorem{lemma}[remark]{Lemma}

\newtheorem*{definition*}{Definition}
\newcommand{\hba}{{_{\beta}H_{\alpha}}}
\newcommand{\kdg}{{_{\delta}K_{\gamma}}}
\newcommand{\lfe}{{_{\phi}L_{\epsilon}}}

\newcommand{\hbak}[1]{{_{\beta_{#1}}H^{#1}_{\alpha_{#1}}}}
\newcommand{\kdgk}[1]{{_{\delta_{#1}}K^{#1}_{\gamma_{#1}}}}

\newcommand{\complex}{\mathbb{C}}

\DeclareMathOperator{\Id}{id}

\DeclareMathOperator{\clspan}{\overline{span}}
\DeclareMathOperator{\Ad}{Ad}

\newcommand{\catC}{\mathbf{C}}
\newcommand{\catD}{\mathbf{D}}

\newcommand{\lt}{\triangleleft}
\newcommand{\rt}{\triangleright}

\newcommand{\hDelta}{\widehat{\Delta}}

\newcommand{\lnspan}{\big[} \newcommand{\rnspan}{\big]}

\newcommand{\frakB}{\mathfrak{B}}
\newcommand{\frakBo}{\frakB^{\dag}}
\newcommand{\frakC}{\mathfrak{C}}

\newcommand{\frakH}{\mathfrak{H}}
\newcommand{\frakK}{\mathfrak{K}}

\newcommand{\cbases}[2]{{_{\mathfrak{#1}}}\mathfrak{#2}_{\mathfrak{#1}^{\dag\!}}} 
\newcommand{\cbaseos}[2]{{_{\mathfrak{#1}^{\dag\!}}}\mathfrak{#2}_{\mathfrak{#1}}} 
\newcommand{\cbasesb}{\cbases{B}{H}}
\newcommand{\cbasesc}{\cbases{C}{K}}

\newcommand{\cbaseosb}{\cbaseos{B}{H}}

\newcommand{\hbeta}{\widehat{\beta}}
\newcommand{\hA}{\widehat{A}}

\newcommand{\mycong}{\xrightarrow{\cong}}
\DeclareMathOperator{\Img}{Im}

\newcommand{\rtensor}[3]{ {_{#1}\! \underset{#2}{\otimes}\! {}_{#3}}}

\newcommand{\htensor}[2]{\rtensor{#1}{\frakH}{#2}}

\newcommand{\rtensorab}{\htensor{\alpha}{\beta}}
\newcommand{\rtensorad}{\htensor{\alpha}{\delta}}

\newcommand{\HadK}{H \rtensorad K}

\newcommand{\rtensorh}{\underset{\frakH}{\otimes}}

\newcommand{\fibre}[3]{ {_{#1}\! \underset{#2}{\ast}\! {}_{#3}}}

\newcommand{\hfibre}[2]{ {_{#1}\! \underset{\frakH}{\ast}\! {}_{#2}}}

\newcommand{\fibreh}{\underset{\frakH}{\ast}}

\newcommand{\fsource}{\htensor{\hbeta}{\alpha}}
\newcommand{\frange}{\rtensorab}

\newcommand{\sfsource}{\fsource}%{\stensor{\hbeta}{\alpha}}
\newcommand{\sfrange}{\frange}%{\stensor{\alpha}{\beta}}

\newcommand{\tl}{\ensuremath \olessthan}
\newcommand{\tr}{\ensuremath \ogreaterthan}

\newcommand{\Hsource}{H \fsource H}
\newcommand{\Hrange}{H \frange H}

\newcommand{\sHsource}{H \sfsource H}
\newcommand{\sHrange}{H \sfrange H}

\newcommand{\HfibreK}{H \rtensorad K}

\newcommand{\Hone}{H \sfsource H \sfsource H}

\newcommand{\Htwo}{H \sfrange H \sfsource  H}

\newcommand{\Hthree}{H \sfrange H \sfrange H}

\newcommand{\Hfour}{(\sHsource) \rtensor{\alpha \lt \alpha}{\frakH}{\beta} H}

\newcommand{\Hfive}{H \rtensor{\hbeta}{\frakH}{\alpha \rt \alpha} (\sHrange)}

\newcommand{\fibreab}{\fibre{\alpha}{\frakH}{\beta}}
\newcommand{\fibread}{\fibre{\alpha}{\frakH}{\delta}}

\newcommand{\AfibreA}{A \fibreab A}

\newcommand{\AfibreB}{A \fibread B}

\newcommand{\kalpha}[1]{|\alpha\rangle_{\leg{#1}}}
\newcommand{\balpha}[1]{\langle\alpha|_{\leg{#1}}}
\newcommand{\kbeta}[1]{|\beta{}\rangle_{\leg{#1}}}
\newcommand{\bbeta}[1]{\langle\beta|_{\leg{#1}}}

\newcommand{\kdelta}[1]{|\delta\rangle_{\leg{#1}}}
\newcommand{\bdelta}[1]{\langle\delta|_{\leg{#1}}}
\newcommand{\khbeta}[1]{|\hbeta{}\rangle_{\leg{#1}}}

\newcommand{\lft}[1]{l(#1)}
\newcommand{\rgt}[1]{r(#1)}
\newcommand{\lfte}[1]{l_{E}(#1)}
\newcommand{\rgte}[1]{r_{E}(#1)}
\newcommand{\lftf}[1]{l_{F}(#1)}
\newcommand{\rgtf}[1]{r_{F}(#1)}
\newcommand{\cfact}{\ensuremath \mathrm{C}^{*}\mathrm{\text{-}fact}}

\DeclareMathOperator{\Dom}{Dom}

\DeclareMathOperator{\paut}{PAut}

\newcommand{\hmgs}[2]{\mathscr{H}_{#2}(#1)}
\newcommand{\hmg}[1]{\mathscr{H}(#1)}
\newcommand{\cL}{{\mathscr{L}}}

\newcommand{\mso}{\mathscr{O}}
\newcommand{\cK}{\mathscr{K}}

\newcommand{\msa}{{\mathscr{A}}} \newcommand{\hmsa}{\mkern
  4mu\widehat{\mkern -4mu\mathscr{A}\mkern 4mu}\mkern -2mu}

\newcommand{\msb}{{\mathscr{B}}}
\newcommand{\msd}{{\mathscr{D}}}
\newcommand{\msc}{{\mathscr{C}}}
\newcommand{\mse}{\mathscr{E}}
\newcommand{\msm}{{\mathscr{M}}}

\newcommand{\msf}{{\mathscr{F}}}
\newcommand{\msg}{{\mathscr{G}}}

\newcommand{\leg}[1]{#1}

\newcommand{\cbalg}{\mathbf{\cbaseosb\text{-}\cbasesb\text{-}alg}}
\newcommand{\bfamily}{(B,\Theta)\mathbf{\text{-}C^{*}\text{-}family}}
\newcommand{\bbimod}{(B,\Theta)\mathbf{\text{-}bimod}}
\newcommand{\rbimod}{B\mathbf{\text{-}bimod}}
\newcommand{\cbbimod}{\mathbf{\cbasesb\text{-} bimod}}

\newcommand{\thetafam}{\mathrm{Fam_{\Theta}}}
\newcommand{\tfam}{\mathrm{Fam}}

\newcommand{\opi}{\tr_{\pi}}

\newcommand{\rE}{{_{r}E}}
\newcommand{\sE}{{_{s}E}}
\newcommand{\Esource}{E {_{s}\tl} E}
\newcommand{\Erange}{E \tr_{r} E}

\newcommand{\funct}[1]{\bm{I}#1}
\newcommand{\functs}[1]{\bm{J}#1}
\DeclareMathOperator{\pinn}{PInn}
\DeclareMathOperator{\pinns}{PInn_{\mathrm{sep}}}
\DeclareMathOperator{\piso}{PIso}

\newcommand{\flip}{\mathrm{flip}}
\title{From Hopf $C^{*}$-families to concrete Hopf $C^{*}$-bimodules}

\author{Thomas Timmermann\\[1ex]
  \texttt{timmermt@math.uni-muenster.de}\\ SFB 478
  ``Geometrische Strukturen in der Mathematik''\\ Hittorfstr.\
  27, 48149 M\"unster}

\date{\today}

\begin{document}
 \xyrequire{matrix} \xyrequire{arrow}

\maketitle 

\abstract{In the setting of von Neumann algebras, measurable quantum
  groupoids have successfully been axiomatized and studied by Enock,
  Vallin, and Lesieur \cite{enock:1,lesieur,vallin:2}, whereas in the setting of $C^{*}$-algebras, a
  similar theory of locally compact quantum groupoids could not yet be
  developed.  Some basic building blocks for such a theory, like
  analogues of a Hopf-von Neumann bimodule and of a
  pseudo-multiplicative unitary, were introduced in
  \cite{timmermann,timmermann:hopf}.  The approach in
  \cite{timmermann,timmermann:hopf}, however, is restricted to
  ``decomposable'' quantum groupoids which generalize $r$-discrete
  groupoids. Recently, we developed a general approach
  \cite{timmer:ckac,timmer:cpmu} that  covers all locally compact
  groupoids. In this article, we explain how the special theory of
  \cite{timmermann,timmermann:hopf} embeds into the general one of
  \cite{timmer:ckac,timmer:cpmu}.  }

\section{Introduction}

In the setting of operator algebras, quantum groupoids have
successfully been axiomatized and studied only on the level of von
Neumann algebras
\cite{enock:9,enock:10,enock:1,enock:2,lesieur,vallin:2} but not on
the finer level of $C^{*}$-algebras.  The theory of locally compact
quantum groups \cite{vaes:1,vaes:30} and the theory of measurable
quantum groupoids \cite{lesieur} suggest that on the level of
$C^{*}$-algebras, a quantum groupoid should be some kind of Hopf
$C^{*}$-bimodule equipped with operator-valued Haar weights, and that
this Hopf $C^{*}$-bimodule is closely related to a
$C^{*}$-pseudo-multiplicative unitary which encodes the
quantum groupoid and its generalized Pontrjagin dual.

A first study of such Hopf $C^{*}$-bimodules and
$C^{*}$-pseudo-multiplicative unitaries was started in 
\cite{timmermann,timmermann:hopf}. But the theory developed there
applies only to a special class of quantum groupoids that are
analogues of $r$-discrete groupoids. Recently, we introduced a general
definition of Hopf $C^{*}$-bimodules and $C^{*}$-pseudo-multiplicative
unitaries \cite{timmer:ckac,timmer:cpmu,timmer:fdcpmu} that, we hope,
should provide the right basis for the study of quantum groupoids on
the level of $C^{*}$-algebras. The purpose of this article is to
explain how the special theory developed in
\cite{timmermann,timmermann:hopf} fits into the general framework
introduced in \cite{timmer:cpmu}.

In both approaches, the definitions of the basic objects involve the
notion of a bimodule over a $C^{*}$-algebra $B$, of a relative tensor
product of bimodules, of (generalized) $C^{*}$-algebras represented on
such bimodules, the fiber product of such (generalized)
$C^{*}$-algebras, and many related constructions. The difference
between the two approaches lies in the choices of the category of
bimodules and the category of represented (generalized)
$C^{*}$-algebras. In \cite{timmermann:hopf}, we have to restrict
ourselves to the special case where the two module structures on the
bimodules and (generalized) $C^{*}$-algebras are related by a family
of partial automorphisms of the underlying $C^{*}$-algebra.

In this article, we construct a functor from the category of bimodules
used in \cite{timmermann:hopf} to the category of bimodules used in
\cite{timmer:cpmu} that is full, faithful, and monoidal in the sense
that it preserves the relative tensor product. Moreover, we construct
a functor from the category of generalized $C^{*}$-algebras (called
$C^{*}$-families) used in \cite{timmermann:hopf} to the category of
concrete $C^{*}$-algebras used in \cite{timmer:ckac,timmer:cpmu} that
is faithful and submonoidal in the sense that it embeds the fiber
product of $C^{*}$-families into the fiber product of concrete
$C^{*}$-algebras.  These two constructions depend on the choice of a
covariant representation of some dynamical system.  Using these 
functors, we associate to suitable Hopf $C^{*}$-families
\cite{timmermann:hopf} concrete Hopf $C^{*}$-bimodules
\cite{timmer:cpmu}, and to suitable pseudo-multiplicative unitaries on
$C^{*}$-modules \cite{timmermann:hopf} $C^{*}$-pseudo-multiplicative
unitaries \cite{timmer:cpmu}.

This work was supported by the SFB 478 ``Geometrische Strukturen in
der Mathematik'' which is funded by the Deutsche
Forschungsgemeinschaft (DFG).

\paragraph{Organization of the article}
This article is organized as follows.

In Section 2, we fix notation and recall some preliminaries concerning
$C^{*}$-modules, partial automorphisms, and monoidal categories.

In Section 3, we  introduce a functor from $C^{*}$-bimodules and
homogeneous operators to Hilbert spaces and ordinary operators which
underlies the constructions in Section 4 and 5.

In Section 4, we  embed the monoidal category
of certain admissible $C^{*}$-bimodules over a $C^{*}$-algebra into
the monoidal category of $C^{*}$-bimodules over a $C^{*}$-base.

In Section 5, we to embed the category of certain admissible
$C^{*}$-families on $C^{*}$-bimodules into the category of
$C^{*}$-algebras over a $C^{*}$-base. This embedding is not
monoidal but embeds the fiber product of $C^{*}$-families into
the fiber product of $C^{*}$-algebras over a $C^{*}$-base.

In Section 6, we use the functor constructed in Section 5 to embed the
category of admissible Hopf $C^{*}$-families over a $C^{*}$-algebra
into the category of concrete Hopf $C^{*}$-bimodules.  Moreover, we
associate to a large class of pseudo-multiplicative unitaries on
$C^{*}$-modules $C^{*}$-pseudo-multiplicative unitaries in such a way
that the legs of these unitaries are related by the functor
constructed in Section 5.

\section{Preliminaries}

Given a subset $Y$ of a normed space $X$, we denote by $[Y] \subseteq
X$ the closed linear span of $Y$.  If $H$ is a Hilbert space and $X
\subseteq {\cal L}(H)$, then $X'$ denotes the commutant of $X$.

We shall make extensive use of (right) $C^{*}$-modules, also known as
Hilbert $C^{*}$-modules or Hilbert modules. A standard reference is
\cite{lance}.

All sesquilinear maps like inner products of Hilbert spaces
or $C^{*}$-modules are assumed to be conjugate-linear in the first
component and linear in the second one.

\paragraph{$C^{*}$-modules}
Let $A$ and $B$ be $C^{*}$-algebras.  
Given $C^{*}$-modules $E$ and $F$ over $B$, we denote the space of all
adjointable operators from $E$ to $F$ by ${\cal L}_{B}(E,F)$.

Let $E$ and $F$ be $C^{*}$-modules over $A$ and $B$, respectively, and
let $\pi \colon A \to {\cal L}_{B}(F)$ be a $*$-homomorphism. Then one
can form the internal tensor product $E \otimes_{\pi} F$, which is a
$C^{*}$-module over $B$ \cite[Chapter 4]{lance}. This $C^{*}$-module
is the closed linear span of elements $\eta \otimes_{A} \xi$, where
$\eta \in E$ and $\xi \in F$ are arbitrary, and $\langle \eta
\otimes_{\pi} \xi|\eta' \otimes_{\pi} \xi'\rangle = \langle
\xi|\pi(\langle\eta|\eta'\rangle)\xi'\rangle$ and $(\eta \otimes_{\pi}
\xi)b=\eta \otimes_{\pi} \xi b$ for all $\eta,\eta' \in E$, $\xi,\xi'
\in F$, and $b \in B$.  We denote the internal tensor product by
``$\tr$''; thus, for example, $E \tr_{\pi} F=E \otimes_{\pi} F$. If
the representation $\pi$ is understood, we write ``$\tr$'' instead of
$"\tr_{\pi}$''.

Given $E$, $F$ and $\pi$ as above, we define a {\em flipped internal
  tensor product} $F {_{\pi}\tl} E$ as follows. We equip the algebraic
tensor product $F \odot E$ with the structure maps $\langle \xi \odot
\eta | \xi' \odot \eta'\rangle := \langle \xi| \pi(\langle
\eta|\eta'\rangle) \xi'\rangle$, $(\xi \odot \eta) b := \xi b \odot
\eta$, and by factoring out the null-space of the semi-norm
$\zeta\mapsto \| \langle \zeta|\zeta\rangle\|^{1/2}$ and taking
completion, we obtain a $C^{*}$-$B$-module $F {_{\pi}\tl} E$.  This is
the closed linear span of elements $\xi {_{\pi}\tl} \eta$, where $\eta \in E$
and $\xi \in F$ are arbitrary, and $\langle \xi {_{\pi}\tl} \eta|\xi'
{_{\pi}\tl} \eta'\rangle = \langle \xi|\pi(\langle\eta|\eta'\rangle)\xi'\rangle$ and
$(\xi {_{\pi}\tl} \eta)b=\xi b {_{\pi}\tl} \eta$ for all $\eta,\eta' \in E$, $\xi,\xi'
\in F$, and $b\in B$. As above, we write ``$\tl$'' instead of
``${_{\pi}\tl}$'' if the representation $\pi$  is understood.

Evidently, the usual and the flipped internal tensor product are related by
a unitary map $\Sigma \colon F \tr E \mycong E \tl F$, $\eta \tr \xi
\mapsto \xi \tl \eta$.

By a {\em right $C^{*}$-$A$-$B$-bimodule} we mean a $C^{*}$-module $E$ 
over $B$ that is full in the sense that $[\langle E|E\rangle]=B$,
together with a fixed nondegenerate representation $A \to {\cal
  L}_{B}(E)$.  Evidently, the class of all right
$C^{*}$-$A$-$B$-bimodules forms a category with respect to the
morphism sets
\begin{align*}
  {\cal L}^{A}_{B}(E,F) = \{ T \in {\cal L}_{B}(E,F) \mid
  Ta\xi=aT\xi \text{ for all } a\in A,\xi \in E\}.
\end{align*}

\paragraph{Partial automorphisms}
Let $B$ be a $C^{*}$-algebra.  A {\em partial automorphism} of $B$ is
a $*$-automorphism $\theta \colon \Dom(\theta) \to \Img(\theta)$,
where $\Dom(\theta)$ and $\Img(\theta)$ are closed ideals of
$B$. Since the composition and the inverse of partial automorphisms
are partial automorphisms again, the set $\paut(B)$ of all partial
automorphisms of $B$ forms an inverse semigroup \cite{paterson}.  For
each $\sigma \in \paut(B)$, we put $\sigma^{*}:=\sigma^{-1}$. Given
$\sigma,\sigma' \in \paut(B)$, we say that $\sigma$ {\em extends}
$\sigma'$ and write $\sigma \geq \sigma'$ if $\Dom(\sigma') \subseteq
\Dom(\sigma)$ and $\sigma|_{\Dom(\sigma')}=\sigma'$.
Given partial automorphisms $\theta,\theta' \in \paut(B)$, we denote
by $\theta \wedge \theta'$ the largest partial automorphism that is
extended by $\theta$ and $\theta'$; thus, $\theta \wedge \theta'
=\theta|_{I}=\theta'|_{I}$, where $I \subseteq \Dom(\theta) \cap
\Dom(\theta')$ denotes the largest ideal on which $\theta$ and
$\theta'$ coincide.

\paragraph{Monoidal categories and functors}

 Let us briefly recall the definition of a monoidal category and of a
monoidal functor; for details, see \cite{maclane}.
A {\em monoidal structure} on a category $\catC$
consists of
\begin{itemize}
\item a bifunctor $\odot \colon \catC \times \catC \to \catC$, $(E,F)
  \mapsto E \odot F$;
\item an object $I \in \catC$ called the unit;
\item for each $E,F,G \in \catC$, an isomorphism
$\alpha_{E,F,G}\colon (E \odot F) \odot G \to E \odot (F
\odot G)$  that is natural in $E,F,G$ and makes the following
diagram commute for all $E,F,G,H \in \catC$:
\begin{align*}
  \xymatrix@R=15pt@C=45pt{ ((E \odot F) \odot G) \odot H
    \ar[r]^{\alpha_{E \odot F,G,H}} \ar[d]^{\alpha_{E,F,G \odot \Id}}
    & (E \odot F) \odot (G \odot H)
    \ar[dd]^{\alpha_{E,F,G\odot H}} \\
    (E \odot (F \odot G) \odot H)
    \ar[d]^{\alpha_{E,F\odot G,H}} & \\
    E \odot ((F \odot G) \odot H) \ar[r]^{\Id\odot \alpha_{F,G,H}} & E
    \odot (F \odot (G \odot H));}
\end{align*}
 \item for each $E \in \catC$,  isomorphisms $l_{E} \colon
  I \odot E \to E$ and $r_{E} \colon E \odot I \to E$ that are
  natural in $E$ and make the following diagram commute for all
  $E,F \in \catC$:
\begin{align*}
 \xymatrix@R=15pt@C=10pt{ (E \odot I) \odot
    F \ar[rr]^{\alpha_{E,I,F}} \ar[rd]_{r_{E}\odot \Id}&& E \odot (I
    \odot F)
    \ar[ld]^{\Id \odot l_{F}} \\
    & E \odot F. &
  }
\end{align*}
\end{itemize}
A {\em monoidal category} is a category equipped with a monoidal
structure.  

A {\em monoidal functor} between monoidal categories
$\catC$ and $\catD$ consists of
\begin{itemize}
\item a functor $\Phi \colon \catC \to \catD$;
\item a natural transformation $\tau \colon \odot \circ (\Phi \times
  \Phi) \to \Phi \circ \odot$  that makes the following diagram commute
  for all $E,F,G \in \catC$:
  \begin{align*}
    \xymatrix@R=15pt@C=45pt{
      (\Phi E \odot \Phi F) \odot \Phi G \ar[r]^{\alpha_{\Phi E,\Phi
          F, \Phi G}}  \ar[d]^{\tau_{E,F} \odot \Id}  &
      \Phi E \odot (\Phi F \odot \Phi G) \ar[d]^{\Id \odot \tau_{F,G}}
      \\
      \Phi(E \odot F) \odot \Phi G 
      \ar[d]^{\tau_{E \odot F, G}}
      & \Phi E \odot \Phi(F \odot G)  \ar[d]^{\tau_{E,F \odot G}}\\
      \Phi((E   \odot F) \odot G) \ar[r]_{\Phi(\alpha_{E,F,G})} &
      \Phi(E \odot (F   \odot G));
    }
  \end{align*}
\item a morphism $\epsilon \colon I_{\catD} \to \Phi(I_{\catC})$,
  where $I_{\catC}$ and $I_{\catD}$ denote the units of $\catC$ and
  $\catD$, respectively,  that makes the following diagrams commute:
  \begin{gather*}
    \xymatrix @R=15pt@C=45pt{
       I_{\catD} \odot \Phi E \ar[r]^{l_{\Phi E}}  \ar[d]^{\epsilon \odot \Id}
       & \Phi E,  \\
       \Phi I_{\catC} \odot \Phi E \ar[r]^{\tau_{I_{\catC},E}} &
       \Phi(I_{\catC} \odot E) \ar[u]^{\Phi(l_{E})}
    }  \qquad
    \xymatrix @R=15pt@C=45pt{
        \Phi E \odot I_{\catD}  \ar[r]^{r_{\Phi E}}  \ar[d]^{\Id \odot \epsilon}
       & \Phi E.  \\
       \Phi E \odot \Phi I_{\catC}  \ar[r]^{\tau_{E,I_{\catC}}} &
       \Phi(E \odot I_{\catC}) \ar[u]^{\Phi(r_{E})}
    }  
  \end{gather*}
\end{itemize}

\section{Untwisting homogeneous operators on $C^{*}$-bimod\-ules}
\label{section:untwist}

The main idea of this article is to use covariant representations of
partial automorphisms to ``untwist'' homogeneous operators and
families of homogeneous operators. We briefly recall the notion of a
homogeneous operator and the appropriate notion of a covariant
representation before we explain the precise construction. 

\paragraph{Homogeneous operators on $C^{*}$-bimodules}
The approach to pseudo-multiplicative unitaries and Hopf
$C^{*}$-bimodules developed in \cite{timmermann,timmermann:hopf} is
based on the concept of right $C^{*}$-bimodules and their operators.
Naturally, this approach leads to operators that do not preserve the
module structure  \cite[Subsection 1.1.1]{timmermann} like,  for example,
\begin{itemize}
\item convolution operators on the $C^{*}$-bimodule of a groupoid;
\item operators of the form $F \to E \tr F$, $\eta' \mapsto
  \xi \tr \eta'$, and $E \to E \tr F$, $\xi' \mapsto
  \xi' \tr \eta$, where $E$ and $F$ are right
  $C^{*}$-bimodules over some $C^{*}$-algebra $B$ and $\xi \in E$,
  $\eta \in F$;
\item the operators comprising a Hopf $C^{*}$-bimodule, if this Hopf
  $C^{*}$-bimodule does not simply correspond to a bundle of quantum
  groups.
\end{itemize}
Therefore, we introduced operators that twist the module structure by
partial automorphisms.  Let us recall the precise definitions.
\begin{definition*}[{\cite{timmermann:hopf}}]
  Let $A,B$ be $C^{*}$-algebras, $E,F$ right
  $C^{*}$-$A$-$B$-bimodules, and $\rho \in \paut(A)$, $\sigma \in
  \paut(B)$.  We call a map $T \colon E \to F$ a {\em
    $(\rho,\sigma)$-homogeneous operator} if
  \begin{enumerate}
  \item $\Img(T) \subseteq [\Img(\rho)F]$ and $Ta\xi=\rho(a)T\xi$ for
    all $a \in \Dom(\rho), \, \xi \in E$, and
  \item there exists a map $S \colon F \to E$ such that $\langle
    SF|E\rangle \subseteq \Dom(\sigma)$ and $\langle\eta|T\xi\rangle =
    \sigma(\langle S\eta|\xi\rangle)$ for all $\xi \in E,\ \eta\in F$.
  \end{enumerate}
  We denote by $\cL^{\rho}_{\sigma}(E,F)$ the set of all
  $(\rho,\sigma)$-homogeneous operators from $E$ to $F$, and put
  $(\cL^{\rho}_{\sigma}(E,F))_{\rho,\sigma}$.
\end{definition*}
Homogeneous operators share many properties of ordinary adjointable
operators on right $C^{*}$-bimodules. Let $A,B$ be $C^{*}$-algebras,
$E,F,G$ right $C^{*}$-$A$-$B$-bimodules, and $\rho,\rho' \in
\paut(A)$, $\sigma,\sigma' \in \paut(B)$.  If $T\colon E \to F$ is a
$(\rho,\sigma)$-homogeneous operator, then $T$ is linear, bounded,
$\|T\|=\|T^{*}\|=\|T^{*}T\|$, and the map $S$ in ii) above is necessarily
unique \cite[Proposition 3.2]{timmermann:hopf}. We call this map $S$
the {\em adjoint} of $T$ and denote it by $T^{*}$. Moreover, by \cite[Proposition 3.7]{timmermann:hopf},
   \begin{align} \label{eq:hmg-op}
  \cL^{\rho'}_{\sigma'}(F,G) \cL^{\rho}_{\sigma}(E,F) &\subseteq
  \cL^{\rho'\rho}_{\sigma'\sigma}(E,G), &
  (\cL^{\rho}_{\sigma}(E,F))^{*} &= \cL^{\rho^{*}}_{\sigma^{*}}(F,E).
\end{align}

We adopt the following notation for families of homogeneous operators.
Let $A,B$ and $E,F,G$ be as above, and let
$\msc=(\msc^{\rho}_{\sigma})_{\rho,\sigma}$ be a family of closed
subspaces $\msc^{\rho}_{\sigma} \subseteq \cL^{\rho}_{\sigma}(E,F)$,
where $\rho \in \paut(B)$, $\sigma \in \paut(A)$.
\begin{itemize}
   \item Given a family $\msd =
     \big(\msd^{\rho}_{\sigma})_{\rho,\sigma}$ of closed subspaces
     $\msd^{\rho}_{\sigma} \subseteq \cL^{\rho}_{\sigma}(E,F)$, we write
     $\msd \subseteq \msc$ if and only if $\msd^{\rho}_{\sigma} \subseteq
     \msc^{\rho}_{\sigma}$ for all $\rho \in \paut(A),\, \sigma \in
     \paut(B)$.
   \item We define a family $\msc^{*} \subseteq \cL(F,E)$ by
     $(\msc^{*})^{\rho}_{\sigma} :=
     \big(\msc^{\rho^{*}}_{\sigma^{*}}\big)^{*}$ for all
     $\rho,\sigma$.
   \item Let $\msd \subseteq \cL(F,G)$ a family of closed subspaces.
     The product $[\msd\msc] \subseteq \cL(E,G)$ is the family given
     by
     \begin{multline*}
       \lnspan\msd\msc\rnspan^{\rho''}_{\sigma''} := \clspan \big\{ T'
       T \, \big|\, T' \in \msd^{\rho'}_{\sigma'},\, T \in
       \msc^{\rho}_{\sigma},\, \rho,\rho' \in \paut(A), \sigma,\sigma'
       \in \paut(B),\\ \rho'\rho \leq \rho'',\, \sigma'\sigma
       \leq \sigma''\big\}
     \end{multline*}
     for all $\rho'' \in \paut(A),\,\sigma'' \in \paut(B)$. Clearly,
     the product $(\msd,\msc)\mapsto [\msd\msc]$ is associative.
     Similarly, we define families $[\msd T],\, [ T'\msc] \subseteq
     \cL(E,G)$ for operators $T \in \cL^{\rho}_{\sigma}(E,F),\, T' \in
     \cL^{\rho'}_{\sigma'}(F,G)$, where $\rho,\rho' \in \paut(A),\,
     \sigma,\sigma' \in \paut(B)$.
   \item We put $[\msc E]:= \clspan \big\{ T \xi \,\big|\, T \in
     \msc^{\rho}_{\sigma}, \rho\in \paut(A), \sigma \in \paut(B), \xi
     \in E\big\}$.
   \item By a slight abuse of notation, we denote by $\msc^{\Id}
     \subseteq \cL(E,F)$ and $\msc_{\Id} \subseteq \cL(E,F)$ the
     families given by
     \begin{align*}
       (\msc^{\Id})^{\rho}_{\sigma} &:= \left\{
         \begin{array}[c]{ll}
           \msc^{\Id}_{\sigma}, & \rho = \Id_{B}, \\
           0, & \text{otherwise},
         \end{array} \right.  &
       (\msc_{\Id})^{\rho}_{\sigma} &:= \left\{
         \begin{array}[c]{ll}
           \msc^{\rho}_{\Id}, & \sigma = \Id_{B}, \\
           0, & \text{otherwise}.
         \end{array} \right. 
     \end{align*}
     Similarly, we define $\cL^{\Id}(E,F)\subseteq \cL(E,F)$ and
     $\cL_{\Id}(E,F)\subseteq \cL(E,F)$.
\end{itemize}

\paragraph{Untwisting homogeneous operators via covariant
  representations}

Till the end of this section, we fix the following data:
\begin{itemize}
\item a $C^{*}$-algebra $B$,
\item an inverse semigroup $\Theta \subseteq \paut(B)$ that satisfies
  $\theta \wedge \theta' \in \Theta$ for all $\theta,\theta' \in \Theta$,
\item a covariant representation $(\pi,\upsilon)$ of $(B,\Theta)$ on a
  Hilbert space $\frakK$, that is, a $*$-homo\-morphism $\pi \colon B
  \to {\cal L}(\frakK)$ and a map $\upsilon \colon \Theta \to {\cal
    L}(\frakK)$ such that 
  \begin{align*}
      \upsilon(\theta)\upsilon(\theta') &= \upsilon(\theta\theta'), &
      \upsilon(\theta)^{*} &= \upsilon(\theta^{*}), &
      \upsilon({\theta})\frakK &=\lnspan
      \pi(\Img(\theta))\frakK\rnspan, &
      \upsilon(\theta)\pi(b)\upsilon(\theta)^{*} &= \pi(\theta(b))
  \end{align*}
  for all $\theta,\theta' \in \Theta$ and $b \in \Dom(\theta)$.  
\end{itemize}
The first two equations above just say that $\upsilon$ is a
homomorphism of inverse semigroups $\Theta \to \piso(\frakK)$, where
$\piso(\frakK)$ denotes the inverse semigroup of all partial
isometries of $\frakK$.

Consider $\frakK$ as a right $C^{*}$-$B$-$\complex$-bimodule via
$\pi$.
  \begin{lemma} \label{lemma:funct-hmg}
    Let $E,F$ be right $C^{*}$-$B$-$B$-bimodules, $\rho,\sigma \in
    \Theta$,  and $T \in
    \cL^{\rho}_{\sigma}(E,F)$.  Then there exists a unique bounded
    linear operator $T \tr_{\pi} \upsilon(\sigma) \colon E \tr_{\pi}
    \frakK \to F \tr_{\pi} \frakK$ such that
    \begin{align*}
  (T \tr_{\pi}
    \upsilon(\sigma))(\xi \tr_{\pi} \zeta) = T\xi \tr_{\pi}
    \upsilon(\sigma)\zeta \quad \text{for all }     \xi \in E, \zeta
    \in \frakK.
    \end{align*}
    If $T \in \cL^{\rho'}_{\sigma'}(E,F) \cap
    \cL^{\rho}_{\sigma}(E,F)$ for some $\rho',\sigma' \in \Theta$,
    then $T \tr_{\pi} \upsilon(\sigma') = T \tr_{\pi}
    \upsilon(\sigma)$.
  \end{lemma}
  \begin{proof}
    By definition, $\upsilon(\sigma) \in \cL^{\sigma}_{\Id}(\frakK)$.
    Hence, the existence of $T \tr_{\pi} \upsilon(\sigma)$ follows
    from \cite[Proposition 5.3]{timmermann:hopf}. Assume that $T \in
    \cL^{\rho'}_{\sigma'}(E,F)$ for some $\rho',\sigma' \in \Theta$.
    Then \cite[Proposition 3.2 (iii)]{timmermann:hopf} implies that $T
    \in \cL^{\rho}_{\sigma''}(E,F)$, where $\sigma''=\sigma \wedge
    \sigma' \in \Theta$.  Let $(u_{\nu})_{\nu}$ be an approximate unit
    of $\Img(\sigma'')$. Since
    $\upsilon(\sigma''\sigma''{}^{*})\frakK=[\pi(\Img(\sigma''))\frakK]$,
    the net $(\pi(u_{\nu}))_{\nu}$ converges in ${\cal L}(\frakK)$
    strongly to $\upsilon(\sigma''\sigma''{}^{*})$. Moreover,
    $T\xi=\lim_{\nu} (T\xi)u_{\nu}$ for each $\xi \in E$
    \cite[Proposition 3.2 (v)]{timmermann:hopf}.  Therefore, we have
    for each $\xi \in E$, $\zeta \in \frakK$,
    \begin{align*}
      T\xi \opi \upsilon(\sigma)\zeta = \lim_{\nu} (T\xi)u_{\nu} \opi
      \upsilon(\sigma)\zeta = \lim_{\nu} T\xi \opi
      \upsilon(\sigma''\sigma''{}^{*}) \upsilon(\sigma)\zeta = T\xi
      \opi \upsilon(\sigma'')\zeta,
    \end{align*}
    whence $T \opi \upsilon(\sigma) = T \opi \upsilon(\sigma'')$. A
    similar argument shows $T \opi \upsilon(\sigma'')=T \opi
    \upsilon(\sigma')$.
  \end{proof}
  \begin{proposition} \label{proposition:funct}
    There exists a functor $\funct{}$ from  the category of right
    $C^{*}$-$B$-$B$-bimodules, where the morphisms are all
    $(\rho,\sigma)$-homogeneous operators with arbitrary $\rho,\sigma \in
    \Theta$, to the category of Hilbert spaces and bounded linear
    operators, given  by
    \begin{itemize}
    \item  $E \mapsto \funct{E} :=  E \tr_{\pi} \frakK$ for each
        right $C^{*}$-$B$-$B$-bimodule $E$, and
      \item $T \mapsto \funct{T} := T \tr_{\pi} \upsilon(\sigma)$ for
        each $(\rho,\sigma)$-homogeneous operator, where $\rho,\sigma
        \in \Theta$.
    \end{itemize}
    This functor satisfies $\funct{(T^{*})}=\funct{(T)}^{*}$ for each
    morphism $T$.
  \end{proposition}
  \begin{proof}
    This follows easily from  Lemma  \ref{lemma:funct-hmg} and
    equation \eqref{eq:hmg-op}.
  \end{proof}

We want to apply the functor above to families of homogeneous
operators that satisfy the following condition:
\begin{definition}
  Let $E,F$ be right $C^{*}$-$B$-$B$-bimodules. A family of closed  subspaces $\msc \subseteq \cL(E,F)$ is {\em $\Theta$-supported} if
  for each $\rho,\sigma \in \paut(B)$, the space
  $\msc^{\rho}_{\sigma}$ is equal to the closed linear span of the
  spaces $\msc^{\rho'}_{\sigma'}$, where $\rho',\sigma' \in \Theta$
  and $\rho' \leq \rho$, $\sigma \leq \sigma'$.
\end{definition}
Given right $C^{*}$-$B$-$B$-bimodules $E$, $F$, denote by $\tfam(E,F)$
the set of all families of closed subspaces $\msc \subseteq \cL(E,F)$,
and by $\thetafam(E,F) \subseteq \tfam(E,F)$ the subset of all
$\Theta$-supported families. For each $\msc \in \tfam(E,F)$, put
\begin{align*}
  \functs{\msc}:= \clspan \{ \funct{c} \mid c \in
  \msc^{\rho}_{\sigma},\rho,\sigma \in \Theta\} \subseteq {\cal
    L}(\funct{E},\funct{F}).
\end{align*}
Inserting the definitions, we find:
\begin{proposition} \label{proposition:functs} Let $E,F,G$ be right
  $C^{*}$-$B$-$B$-bimodules. Then
    \begin{align*}
      \functs{(\msc^{*})} &= \functs{(\msc)}^{*}, &
      [\functs{(\msd)}\functs{(\msc})] &\subseteq \functs{[\msd\msc]}
    \end{align*}
     for all $\msc \in \tfam(E,F)$, $\msd \in \tfam(F,G)$, and
     \begin{align*}
      \msc^{*} &\in \thetafam(F,E), & [\msd\msc]&\in \thetafam(E,G), &
[\functs{(\msd)}\functs{(\msc})] &=      \functs{[\msd\msc]}
    \end{align*}
    for all $\msc \in \thetafam(E,F)$, $\msd \in \thetafam(F,G)$.
    \qed
\end{proposition}

\section{From $C^{*}$-bimodules over $C^{*}$-algebras to
  $C^{*}$-bimod\-ules over $C^{*}$-bases}

In this section, we use the assignments $\funct{}$ and $\functs{}$ to
construct a functor from the category of right $C^{*}$-bimodules over
$C^{*}$-algebras used in \cite{timmermann,timmermann:hopf} to the
category of $C^{*}$-bimodules over $C^{*}$-bases used in
\cite{timmer:ckac,timmer:cpmu}, and show that this functor preserves
the relative tensor product.

\subsection{The monoidal category of $\Theta$-admissible
  $C^{*}$-bimodules over a $C^{*}$-algebra}

The theory developed in \cite{timmermann,timmermann:hopf} is based on
the notion of decomposable right $C^{*}$-bimodules over a
$C^{*}$-algebra and on the internal tensor product of right
$C^{*}$-bimodules. In this subsection, we review the category of
decomposable right $C^{*}$-bimodules, introduce the subcategory of
$\Theta$-admissible right $C^{*}$-bimodules, and explain the monoidal
structure of these categories that is induced by the internal tensor
product.

Throughout this subsection, we fix a $C^{*}$-algebra $B$.

\paragraph{Homogeneous elements of $C^{*}$-bimodules}

\begin{definition*}[{\cite{timmermann:hopf}}]
  Let $E$ be a right $C^{*}$-$B$-$B$-bimodule and $\theta \in
  \paut(B)$.  We call an element $\xi \in E$ {\em
    $\theta$-homogeneous} if $\xi \in E\Dom(\theta)$ and $\xi
  b=\theta(b)\xi$ for all $b \in \Dom(\theta)$.  We denote by
  $\hmgs{E}{\theta}$ the space of all $\theta$-homogeneous elements of
  $E$, and call $E$ {\em decomposable} if the family
  $\hmg{E}:=(\hmgs{E}{\theta'})_{\theta' \in \paut(B)}$ is linearly
  dense in $E$.
\end{definition*}
We also consider $B$ as a right $C^{*}$-$B$-$B$-bimodule and denote by
$\hmgs{B}{\theta} \subseteq B$ the subspace of $\theta$-homogeneous
elements.  By \cite[Proposition 3.14]{timmermann:hopf}, we have for
each right $C^{*}$-$B$-$B$-bimodule $E,F$ and each $\rho,\sigma,\theta
\in \paut(B)$:
  \begin{gather} \label{eq:hmg}
    \begin{gathered}
      \begin{aligned}
        \langle \hmgs{E}{\theta}|\hmgs{E}{\theta'}\rangle &\subseteq
        \hmgs{B}{(\theta^{*}\theta')}, &
        \hmgs{B}{\rho}\hmgs{E}{\theta}\hmgs{B}{\sigma} &\subseteq
        \hmgs{E}{(\rho\theta\sigma)},
      \end{aligned}\\
      \hmgs{E}{\theta} \hmgs{B}{(\theta^{*}\theta)} = \hmgs{E}{\theta}
      =
      \hmgs{B}{(\theta\theta^{*})}\hmgs{E}{\theta}, \\
      \begin{aligned}
        \cL^{\rho}_{\sigma}(E,F)\hmgs{E}{\theta} &\subseteq
        \hmgs{F}{(\rho\theta\sigma^{*})}, & \theta(\hmgs{B}{\theta'}
        \cap \Dom(\theta)) &\subseteq
        \hmgs{B}{(\theta\theta'\theta^{*})}.
      \end{aligned}
    \end{gathered}
  \end{gather}
 \begin{definition}
   Let $\Theta \subseteq \paut(B)$ be an inverse semigroup. We call a
   right $C^{*}$-$B$-$B$-bimodule $E$
   \begin{itemize}
   \item {\em $\Theta$-decomposable} if the family of subspaces
     $\hmgs{E}{\theta}$, where $\theta \in \Theta$, is linearly dense in
     $E$;
   \item {\em $\Theta$-supported} if for each $\theta \in \paut(B)$,
     the space $\hmgs{E}{\theta}$ is the closed linear span of all
     subspaces $\hmgs{E}{\theta'}$, where $\theta' \in \Theta$ and
     $\theta' \leq \theta$.
   \end{itemize}
  \end{definition}
  The two conditions are related to each other as follows:
  \begin{proposition} \label{proposition:supported} Let $\Theta
    \subseteq \paut(B)$ be an inverse semigroup such that $B$ is
    $\Theta$-supported, and let $E$ be a right
    $C^{*}$-$B$-$B$-bimodule. Then $E$ is $\Theta$-decomposable if and
    only if it is decomposable and $\Theta$-supported.
  \end{proposition}
  \begin{proof}
   If $E$ is decomposable and $\Theta$-supported, then clearly
    $E$ is $\Theta$-decomposable. Conversely, assume that $E$ is
    $\Theta$-decomposable. Then $E$ is decomposable, and we have to
    show that it is $\Theta$-supported. For each $\theta \in
    \paut(B)$, denote by $E^{0}_{\theta} \subseteq \hmgs{E}{\theta}$ the
    closed linear span of the spaces $\hmgs{E}{\theta'}$, where $\theta'
    \in \Theta$ and $\theta' \leq \theta$. Then the family
    $(E_{\theta})_{\theta}$ is linearly dense in $E$, and since $B$ is
    $\Theta$-supported, $E_{\theta}^{0}\hmgs{B}{\sigma} \subseteq
    E_{\theta\sigma}^{0}$ for all $\theta,\sigma \in \paut(B)$. Now,
    \cite[Proposition 3.15]{timmermann:hopf} implies
    $E^{0}_{\theta}=\hmgs{E}{\theta}$ for all $\theta \in \paut(B)$.
\end{proof}

Let us describe a natural inverse semigroup $\Theta \subseteq
\paut(B)$ for which $B$ is $\Theta$-supported.
\begin{definition}
  A partial automorphism $\theta$ of a $C^{*}$-algebra $B$ is {\em
(separable)    inner } if there exist a (separable) ideal $I \subseteq
  Z(B)$ and a unitary $u \in M(IB)$ such that
  $\Dom(\theta)=IB=\Img(\theta)$ and $\theta(b)=ubu^{*}$ for all $b
  \in IB$. In that case, we write $\theta=\Ad_{u}$. We denote by
  $\pinn_{(\mathrm{sep})}(B) \subseteq \paut(B)$ the subset of all
(separable)  inner partial automorphisms of $B$.
\end{definition}
\begin{lemma}
  $\pinn(B)$ and $\pinns(B)$ are inverse
  subsemigroups of $\paut(B)$.
\end{lemma}
\begin{proof}
  Let $I,J \subseteq Z(B)$ be ideals and let $u \in M(IB), v \in
  M(JB)$ be unitaries. Then $(\Ad_{u})^{*} = \Ad_{u^{*}} \in
  \pinn(B)$. Since $I\cap J \subseteq Z(IB \cap JB)$, the unitaries
  $u,v$ restrict to unitaries $u',v' \in M(IB \cap JB)$. Now, $\Ad_{u}
  \circ \Ad_{v} = \Ad_{u'v'} \in \pinn(B)$. If $I$ and $J$ are
  separable, then also $I\cap J$ is separable.
\end{proof}
  \begin{proposition} \label{proposition:b-supported}
    The right $C^{*}$-$B$-$B$-bimodule $B$ is $\pinns(B)$-supported.
  \end{proposition}
  \begin{proof}
    By \cite[Proposition 3.19 (iii)]{timmermann:hopf}, there exists
    for each $\theta \in \paut(B)$ and $b \in \hmgs{B}{\theta}$ a
    $\theta' \in \pinns(B)$ such that $\theta' \leq \theta$ and $b \in
    \hmgs{B}{\theta'}$. 
  \end{proof}

The definition of $\Theta$-admissible $C^{*}$-bimodules and the
construction of the monoidal functor involve several ket-bra operators
associated to homogeneous elements.  Let $E$ be a right
$C^{*}$-$B$-$B$-bimodule,   $\theta \in \paut(B)$, and $\xi \in
\hmgs{E}{\theta}$, $b \in \hmgs{B}{\theta}$. By \cite[Propositions 3.12,
3.21]{timmermann:hopf}, there exist homogeneous operators
\begin{align*}
 \lft{\xi}=|\xi\rangle &\in \cL^{\theta}_{\Id}(B,E), &
\rgt{\xi}=|\xi] &\in \cL_{\theta^{*}}^{\Id}(B,E), &
\lfte{b} &\in \cL^{\theta}_{\Id}(E), &  \rgte{b} &\in
  \cL^{\Id}_{\theta^{*}}(E)
\end{align*}
such that for all $b' \in B$ and $\xi' \in E$,
\begin{align*}
  \lft{\xi} b'&= \xi b', & \rgt{\xi} b' &= b'\xi, & \lfte{b} \xi'&=
  b\xi', &
  \rgte{b}\xi' &= \xi' b, \\
  \lft{\xi}^{*}\xi' &= \langle\xi|\xi'\rangle, & \rgt{\xi}^{*}\xi' &=
  \theta(\langle\xi|\xi'\rangle), & \lfte{b}^{*} \xi' &= b^{*}\xi', &
  \rgte{b}^{*}\xi' &= \xi' b^{*}.
  \end{align*}
  We define families
  \begin{align} \label{eq:lr-families}
      \lft{\hmg{E}} &\subseteq \cL_{\Id}(B,E), & \rgt{\hmg{E}}
      &\subseteq \cL^{\Id}(B,E), 
\end{align}
such that $\lft{\hmg{E}}^{\theta}_{\Id} = \lft{\hmgs{E}{\theta}}$ and
$\rgt{\hmg{E}}^{\Id}_{\theta} = \rgt{\hmgs{E}{\theta^{*}}}$ for each
$\theta \in \paut(B)$. Similarly, we define families $ \lfte{\hmg{B}} \subseteq
\cL_{\Id}(E)$ and $\rgte{\hmg{B}} \subseteq \cL^{\Id}(E)$.
\begin{lemma} \label{lemma:lr-module} Let $E$ be a right
  $C^{*}$-$B$-$B$-bimodule. Then
    \begin{align*} 
      \begin{aligned}
        \rgt{\xi}\rgt{b} &= \rgt{b\xi}, &
        \rgte{b}\rgt{\xi} &=  \rgt{\xi b}, 
        & \rgt{\xi} \lft{b} &=
        \lfte{b}\rgt{\xi}, & \rgt{\xi}^{*}\rgt{\xi'} &=
        \rgt{\rgt{\xi}^{*}\xi'}, \\
        \lft{\xi}\lft{b} &= \lft{\xi b}, &
        \lfte{b}\lft{\xi} &= \lft{b\xi}, 
        & \lft{\xi}\rgt{b} &=
        \rgte{b}\lft{\xi}, & \lft{\xi}^{*}\lft{\xi'} &=
        \lft{\lft{\xi}^{*}\xi'}, \\
      && S\lft{\xi} &= \lft{S\xi}, &
      T\rgt{\xi} &= \rgt{T\xi},
      \end{aligned}
    \end{align*}
    for all homogeneous $b \in B, \xi,\xi' \in E$ and $S \in
    \cL^{\theta}_{\Id}(E)$, $T \in \cL^{\Id}_{\theta}(E)$, $\theta \in
    \paut(B)$. Put $\mse:=\hmg{E}$ and $\msb:=\hmg{B}$. If $B$ is
    decomposable, then
    \begin{gather} \label{eq:lr-module}
     \begin{aligned}
 {}       [\rgt{\mse}\rgt{\msb}] &= \rgt{\mse} =
        [\rgte{\msb}\rgt{\mse}], &
        [\rgt{\mse}^{*}\rgt{\mse}] &\subseteq \rgt{\msb},  &
        [\cL^{\Id}(E)\rgt{\mse}] &= \rgt{\mse}, \\
        [\lft{\mse}\lft{\msb}] &= \lft{\mse} =
        [\lfte{\msb}\lft{\mse}], & [\lft{\mse}^{*}\lft{\mse}]
        &\subseteq \lft{\msb}, &
        [\cL_{\Id}(E)\lft{\mse}] &= \lft{\mse}.
      \end{aligned}
    \end{gather}
  \end{lemma}
  \begin{proof}
    The first set of equations can be verified by straightforward
    calculations; we only prove $\rgt{\xi}^{*}\rgt{\xi'} =
    \rgt{\theta(\langle\xi|\xi'\rangle)}$: for all $b \in B$,
    \begin{align*}
      \rgt{\xi}^{*}\rgt{\xi'}b = \rgt{\xi}^{*}(b\xi') =
      (\lfte{b^{*}}\rgt{\xi})^{*}\xi' = (\rgt{\xi}\lft{b^{*}})^{*}\xi'
      = b (\rgt{\xi}^{*}\xi') = \rgt{\rgt{\xi}^{*}\xi'}b.
    \end{align*}
 The equations in \eqref{eq:lr-module} follow
    directly from the equations above, equation \eqref{eq:hmg}, and
    the fact that $Z(B)=\hmgs{B}{\Id}\subseteq B$ is nondegenerate
    \cite[Proposition 3.20 (v)]{timmermann:hopf}.
  \end{proof}

  To construct the monoidal functor, we shall apply the map
  $\functs{}$ constructed in Section 2 to the families in
  \eqref{eq:lr-families}. Then, the following condition on the inverse
  semigroup $\Theta$ turns out to be useful: 
  \begin{definition}
    We call an inverse semigroup $\Theta \subseteq \paut(B)$ {\em
      admissible} if $B$ is $\Theta$-supported, $\Id_{B} \in \Theta$,
    and $\theta \wedge \theta' \in \Theta$ for all $\theta,\theta' \in
    \Theta$.
  \end{definition}
  \begin{definition}
    Let $\Theta \subseteq \paut(B)$ be an admissible inverse
    semigroup.  A right $C^{*}$-$B$-$B$-bimodule $E$ is {\em
      $\Theta$-admissible} if it is $\Theta$-decomposable and
    \begin{align*}
      [\lft{\hmg{E}}^{*}\lft{\hmg{E}}] &= \lft{\hmg{B}}, &
      [\rgt{\hmg{E}}^{*}\rgt{\hmg{E}}] &= \rgt{\hmg{B}}.
    \end{align*}
\end{definition}
Proposition \ref{proposition:supported} immediately  implies:
\begin{corollary} \label{corollary:lr-supported} Let $\Theta \subseteq
  \paut(B)$ be an admissible inverse semigroup and let $E$ be a
  $\Theta$-admissible right $C^{*}$-$B$-$B$-bimodule.  Then the
  families in \eqref{eq:lr-families} are $\Theta$-supported.  \qed
\end{corollary}

\paragraph{The internal tensor product}
The category $\rbimod$ of all right $C^{*}$-$B$-$B$-bimodules carries
a monoidal structure, where $B$ is the unit and for all objects
$E,F,G$ and all morphisms $S,T$,
\begin{itemize}
\item $E \odot F := E \tr  F$ and $S \odot T := S \tr T$,
\item $\alpha_{E,F,G} \colon (E \tr F) \tr G \to E
  \tr (F \tr G)$ is given by $(\eta \tr \xi)
  \tr \zeta \mapsto \eta \tr (\xi \tr \zeta)$,
\item $l_{E}\colon B \tr E \to E$, $r_{E} \colon E \tr B \to E$ are
  given by $b \tr \xi \mapsto b\xi$, $\xi \tr b \mapsto \xi b$,
  respectively.
\end{itemize}

Let $\Theta \subseteq \paut(B)$ be an admissible inverse semigroup. We
shall show that the internal tensor product of $\Theta$-admissible
right $C^{*}$-$B$-$B$-bimodules $E$ and $F$ is $\Theta$-admissible
again.  The proof involves the following generalized ket-bra
operators.  By \cite[Proposition 3.13]{timmermann:hopf}, there exist
for each $\theta \in \paut(B)$, $\xi \in \hmgs{E}{\theta}$, $\eta \in
\hmgs{F}{\theta}$ operators
\begin{align*}
  |\xi\rangle_{1} \in \cL^{\theta}_{\Id}(F,E \tr F),   & \
  |\xi\rangle_{2}  \in \cL^{\theta}_{\Id}(F,   F \tl E), &
  |\eta]_{2}   \in \cL^{\Id}_{\theta^{*}}(E, E \tr F), & \
  |\eta]_{1}  \in \cL^{\Id}_{\theta^{*}}(E, F \tl E)
\end{align*}
such that for all $\xi' \in E$ and $\eta' \in F$
\begin{gather*}
  \begin{aligned}
    |\xi\rangle_{1} \eta'&= \xi \tr \eta', & |\xi\rangle_{2} \eta' &=
    \eta' \tl \xi, & |\eta]_{2} \xi' &= \xi' \tr \eta, & |\eta]_{1}
    \xi' &= \eta \tl \xi',
  \end{aligned}
  \\
  \begin{aligned}
    \langle \xi|_{1}(\xi' \tr \eta')&=\langle\xi|\xi'\rangle\eta'=\langle
    \xi|_{2}(\eta' \tl \xi'), &
    [\eta|_{2}(\xi' \tr \eta') &= \xi'
    \theta(\langle\eta|\eta'\rangle) = [\eta|_{1}(\eta' \tl \xi').
  \end{aligned}
\end{gather*}
We define families
\begin{align} \label{eq:ketbra-families}
  |\hmg{E}\rangle_{1} &\subseteq \cL_{\Id}(F,E \tr F), &
|\hmg{F}]_{2} &\subseteq \cL^{\Id}(E,E \tr F)
\end{align}
such that
$(|\hmg{E}\rangle_{1})^{\theta}_{\Id}=|\hmgs{E}{\theta}\rangle_{1}$
and $(|\hmg{F}]_{2})^{\Id}_{\theta} = |\hmgs{F}{\theta^{*}}]_{2}$ for
all $\theta \in \paut(B)$.
\begin{proposition} \label{proposition:itp} Let $E$ and $F$ be
  $\Theta$-admissible right $C^{*}$-$B$-$B$-bimodules. Then also the
  right $C^{*}$-$B$-$B$-bimodule $E \tr F$ is $\Theta$-admissible.
  For each $\theta'' \in \paut(B)$, the space $\hmgs{E \tr
    F}{\theta''}$ is the closed linear span of the subspaces
  $\hmgs{E}{\theta} \tr \hmgs{F}{\theta'}$, where $\theta,\theta' \in
  \Theta$ and $\theta\theta' \leq \theta''$, in particular,
  \begin{align} \label{eq:hmg-2}
    \lft{\hmg{E \tr F}}&=[|\hmg{E}\rangle_{1}\lft{\hmg{F}}], &
    \rgt{\hmg{E \tr F}} &= [|\hmg{F}]_{2}\rgt{\hmg{E}}]. 
  \end{align}
\end{proposition}
\begin{proof}
  By \cite[Proposition 3.17]{timmermann:hopf}, the space $\hmgs{E \tr
    F}{\theta''}$ is for each $\theta'' \in \paut(B)$ the closed
  linear span of the subspaces $\hmgs{E}{\theta} \tr
  \hmgs{F}{\theta'}$, where $\theta,\theta' \in \paut(B)$ and
  $\theta\theta' \leq \theta''$. Equation \eqref{eq:hmg-2} follows.
  Since $E$ and $F$ are $\Theta$-supported, the same statement holds
  if we allow $\theta,\theta'$ to take values in $\Theta$ only.  Put
  $\mse:=\hmg{E}$, $\msf:=\hmg{F}$, $\msg:=\hmg{E \tr F}$, and
  $\msb:=\hmg{B}$. 
Using Lemma
  \ref{lemma:lr-module} and the assumptions on $E$ and $F$, we find
  \begin{align*}
    [\lft{\msg}^{*}\lft{\msg}] = [\lft{\msf}^{*}\langle\mse|_{1}
    |\mse\rangle_{1}\lft{\msf}] &=
    [\lft{\msf}^{*}\lfte{\msb}\lft{\msf}] \\
    &= [\lft{\msf}^{*}\lft{\msb\msf}] = [\lft{\msf}^{*}\lft{\msf}]
    =
    \lft{\msb}, \\
    [\rgt{\msg}^{*}\rgt{\msg}] =
    [\rgt{\mse}^{*}[\msf|_{2}|\msf]_{2}\rgt{\mse}]
    &= [\rgt{\mse}^{*} \rgte{\msb}\rgt{\mse}] \\
    &= [\rgt{\mse}^{*}\rgt{[\mse\msb]}]= [\rgt{\mse}^{*}\rgt{\mse}]
    = \rgt{\msb}. \qedhere
  \end{align*}
\end{proof}
\begin{remark}
  Given a right $C^{*}$-$B$-$B$-bimodule $G$ and homogeneous elements
  $\zeta,\zeta' \in G$, define $ [\zeta|\zeta'] :=
  \rgt{\zeta}^{*}\zeta'\in B$.  Then for all right
  $C^{*}$-$B$-$B$-bimodules $E,F$ and all homogeneous $\xi,\xi'\in E$,
  $\eta,\eta' \in F$, we have $[\xi \tr \eta |\xi' \tr \eta'] =
  \rgt{\xi}^{*} [\eta|_{2}(\xi' \tr \eta') = \rgt{\xi}^{*}(\xi'
  \rgt{\eta}^{*}\eta') = [\xi|\xi' [\eta|\eta']]$.  This formula
  resembles the formula $\langle \xi \tr \eta|\xi' \tr \eta'\rangle =
  \langle \xi| \langle \eta|\eta'\rangle\xi'\rangle$ used in the
  definition of the internal tensor product of right $C^{*}$-bimodules
  and is the natural choice for the definition of the internal tensor
  product of {\em left} $C^{*}$-bimodules.
  \end{remark}

The main result of this subsection is:
\begin{corollary}\label{corollary:bbimod-monoidal}
  Let $\Theta \subseteq \paut(B)$ be an admissible inverse semigroup.
  Then the full subcategory $\bbimod$ of $\rbimod$ that consists of
  all $\Theta$-admissible right $C^{*}$-$B$-$B$-bimodules is monoidal.
  \qed
\end{corollary}

\begin{remark} \label{remark:bbimod-flip}
  The flipped internal tensor product  defines another monoidal
  structure on the category $\bbimod$, where $B$ is the unit again and for all
  objects $E,F,G$ and all morphisms $S,T$,
  \begin{itemize}
  \item $E \odot^{\flip} F := E \tl  F$ and $S \odot^{\flip} T := S \tl T$,
  \item $\alpha^{\flip}_{E,F,G} \colon (E \tl F) \tl G \to E
    \tl (F \tl G)$ is given by $(\eta \tl \xi)
    \tl \zeta \mapsto \eta \tl (\xi \tl \zeta)$,
  \item $l^{\flip}_{E}\colon B \tl E \to E$, $r^{\flip}_{E} \colon E \tl B \to E$ are
    given by $b \tl \xi \mapsto \xi b$, $\xi \tl b \mapsto b\xi$,
    respectively.
  \end{itemize}
  Denote this monoidal category by $\bbimod^{\flip}$. 
\end{remark}
\subsection{The monoidal category of $C^{*}$-bimodules over a
    $C^{*}$-base}

  The theory developed in \cite{timmer:ckac,timmer:cpmu} is based on
  the notion of a $C^{*}$-bimodule over a $C^{*}$-base, which in turn
  is based on the notion of a $C^{*}$-factorization.  
  \begin{definition*}[{\cite{timmer:cpmu}}]
    A {\em $C^{*}$-base} is a triple $(\frakH,\frakB,\frakBo)$,
    shortly written $\cbases{B}{H}$, consisting of a Hilbert space
    $\frakH$ and two commuting nondegenerate $C^{*}$-algebras
    $\frakB,\frakB^{\dag} \subseteq {\cal L}(\frakH)$.

    Let $H$ be a Hilbert space and $\cbasesb$ a $C^{*}$-base.  A {\em
      $C^{*}$-factorization} of $H$ with respect to $\cbases{B}{H}$ is
    a closed subspace $\alpha \subseteq {\cal L}(\frakH,H)$ satisfying
    $[\alpha^{*}\alpha] = \frakB$, $[\alpha \frakB] = \alpha$, and
    $[\alpha \frakH] = H$. We denote the set of all such
    $C^{*}$-factorizations by $\cfact(H;\cbases{B}{H})$.
  \end{definition*}

  Let $\alpha$ be a $C^{*}$-factorization of a Hilbert space $H$ with
  respect to a $C^{*}$-base $\cbases{B}{H}$.  Then $\alpha$ is a
  concrete $C^{*}$-module and a full right $C^{*}$-module over
  $\frakB$ with respect to the inner product $\langle \xi|\xi'\rangle
  :=\xi^{*}\xi'$.  Moreover, we can identify
  $\alpha \tr \frakH$ and $\frakH \tl \alpha$ with $H$ via the
  unitaries
\begin{align} \label{eq:rtp-iso} m_{\alpha} \colon \alpha \tr \frakH
  &\to H, \ \xi \tr \zeta \mapsto \xi\zeta, & m_{\alpha}^{op} \colon
  \frakH \tl \alpha &\to H, \ \zeta \tl \xi \mapsto \xi\zeta,
\end{align}
and there exists a nondegenerate and faithful representation
$\rho_{\alpha} \colon \frakB^{\dag} \to  {\cal L}(H)$ such that
\begin{align*}
  \rho_{\alpha}(b^{\dag})\xi\zeta = \xi b^{\dag}\zeta \quad \text{for
    all } b^{\dag} \in \frakB^{\dag}, \, \xi\in \alpha, \, \zeta \in
  \frakH;
\end{align*}
see \cite[Subsection 2.1]{timmer:cpmu}. 

Let $\beta$ be a $C^{*}$-factorization of a Hilbert space  $K$ with
respect to $\cbasesb$. We put
\begin{align*}
  {\cal L}(H_{\alpha},K_{\beta}) := \{ T \in {\cal L}(H,K) \mid
  T\alpha \subseteq \beta, \, T^{*}\beta \subseteq \alpha\}.
\end{align*}
Each $T \in {\cal L}(H_{\alpha},K_{\beta})$ defines an
operator $T_{\alpha} \in {\cal L}_{\frakB}(\alpha,\beta)$ by $\xi
\mapsto T\xi$ with adjoint $\eta \mapsto T^{*}\eta$. Moreover, the relation
$T\rho_{\alpha}(b^{\dag})\xi\zeta = T\xi b^{\dag} \zeta = 
\rho_{\beta}(b^{\dag})T \xi\zeta$, valid  for all $\xi \in \alpha, \zeta \in
\frakH$ implies
\begin{align} \label{eq:trho-commute}
  T\rho_{\alpha}(b^{\dag})=\rho_{\beta}(b^{\dag})T \quad \text{for all
  } T \in {\cal L}(H_{\alpha},K_{\beta}),  b^{\dag} \in \frakBo.
\end{align}

Let $\cbasesc$ be a $C^{*}$-base.  We call a $C^{*}$-factorization
$\beta \in \cfact(H;\cbasesc)$ {\em compatible} with $\alpha$, written
$\alpha \perp \beta$, if $\lnspan \rho_{\alpha}(\frakBo)\beta\rnspan
=\beta $ and $\lnspan
\rho_{\beta}(\frakC^{\dag})\alpha\rnspan=\alpha$. In that case,
$\rho_{\alpha}(\frakBo) \subseteq {\cal L}(H_{\beta})$ and
$\rho_{\beta}(\frakC^{\dag}) \subseteq {\cal L}(H_{\alpha})$; in
particular, these $C^{*}$-algebras commute.

\begin{definition}
  Let $\cbasesb$ and $\cbasesc$ be $C^{*}$-bases. A {\em
    $C^{*}$-$\cbasesb$-$\cbasesc$-bimodule} is a triple $(H,\beta,\alpha)$,
  briefly denoted by ${_{\beta}H_{\alpha}}$, consisting of a Hilbert
  space $H$ and compatible $C^{*}$-factorizations $\alpha \in
  \cfact(H;\cbasesc)$ and $\beta \in \cfact(H;\cbaseosb)$. 
\end{definition}

Let $\cbasesb$ and $\cbasesc$ be $C^{*}$-bases as before. Moreover,
let $H$ and $K$ be Hilbert spaces with $C^{*}$-factorizations $\alpha
\in \cfact(H;\cbasesb)$ and $\delta \in \cfact(K;\cbaseosb)$. The {\em
  $C^{*}$-relative tensor product} of $H$ and $K$ with respect to
$\alpha$ and $\delta$ is the Hilbert space
\begin{align*}
  \HadK := \alpha \tr \cbasesb \tl \delta.
\end{align*}
 The unitaries \eqref{eq:rtp-iso} induce isomorphisms 
  \begin{gather} \label{eq:rtp-space}
    \begin{gathered}
      \alpha \tr_{\rho_{\delta}} K \xleftarrow{\Id \tr
        m_{\delta}^{op}} \HadK
      \xrightarrow{m_{\alpha} \tl \Id} H {_{\rho_{\alpha}}\tl} \beta, \\
      \xi \tr_{\rho_{\delta}} \eta\zeta \mapsfrom \xi \tr_{\frakB} \zeta
 \,      {_{\frakBo}\tl} \eta \mapsto \xi\zeta {_{\rho_{\alpha}} \tl} \eta.
    \end{gathered}
  \end{gather}
Using these isomorphisms, we define for each $\xi \in \alpha$ and
$\eta \in \delta$ two pairs of adjoint operators
\begin{align*}
  |\xi\rangle_{\leg{1}} \colon K &\to \HfibreK, \ \zeta \mapsto \xi
  \tr \omega, & \langle \xi|_{\leg{1}}:=|\xi\rangle_{\leg{1}}^{*}\colon
  \xi' \tr \omega &\mapsto
  \rho_{\delta}(\langle\xi|\xi'\rangle)\omega, \\
  |\eta\rangle_{\leg{2}} \colon H &\to \HfibreK, \ \omega \mapsto \omega
  \tl \eta, & \langle\eta|_{\leg{2}} := |\eta\rangle_{\leg{2}}^{*}
  \colon \omega \tl\eta &\mapsto \rho_{\alpha}(\langle
  \eta|\eta'\rangle)\omega.
\end{align*}
We put $\kalpha{1} := \big\{ |\xi\rangle_{\leg{1}} \,\big|\, \xi \in
\alpha\big\}$ and similarly define $\balpha{1}$, $\kdelta{2}$,
$\bdelta{2}$.

For each $C^{*}$-factorization $\beta \in \cfact(H;\cbasesc)$ that is
compatible with $\alpha$, the space
\begin{align*}
  \beta \lt \delta := [|\delta\rangle_{2}\beta] \subseteq {\cal
    L}(\frakK,\HadK)
\end{align*}
is a $C^{*}$-factorization of $\HadK$ with respect to $\cbasesc$; likewise, 
for each  $C^{*}$-factorization $\gamma \in
\cfact(H;\cbasesc)$ that is compatible with $\gamma$, the space
\begin{align*}
  \alpha \rt \gamma := [|\alpha\rangle_{1}\gamma] \subseteq {\cal
    L}(\frakK,\HadK)
\end{align*}
is a $C^{*}$-factorization of $\HadK$ with respect to $\cbasesc$
\cite[Proposition 2.7]{timmer:cpmu}.

Let $L$ and $M$ be Hilbert spaces with $C^{*}$-factorizations $\beta
\in \cfact(L;\cbasesb)$ and $\gamma \in \cfact(M;\cbaseosb)$. Then
there exists for each $S \in {\cal L}(H_{\alpha},L_{\beta})$, $T
\in {\cal L}(K_{\delta},M_{\gamma})$ an operator
\begin{align*}
  S \rtensorh T := S_{\alpha} \tr \Id \tl T_{\delta} \in {\cal L}(H
  \htensor{\alpha}{\delta} K, L \htensor{\beta}{\gamma} M).
\end{align*}

The class of all $C^{*}$-$\cbasesb$-$\cbasesc$-bimodules forms a
category with respect to the morphisms
\begin{align*}
  {\cal L}({_{\beta}H_{\alpha}},{_{\delta}K_{\gamma}}) := {\cal
    L}(H_{\alpha},K_{\gamma}) \cap {\cal L}(H_{\beta},K_{\delta}),
\end{align*}
as one can easily verify.  If $\cbasesc=\cbasesb$, we denote this
category by $\cbbimod$.

\begin{theorem} \label{theorem:cbbimod}
  The category $\cbbimod$ carries a structure of a monoidal category, where
  $\cbaseosb$ is the unit and for all objects $\hba, \kdg,\lfe$ and all
  morphisms $S,T$,
  \begin{itemize}
  \item $\hba \odot \kdg = \hba \rtensorh \kdg := {_{\beta \lt \delta}(H
      \htensor{\alpha}{\delta} K)_{\alpha \rt \gamma}}$ and $S \odot T
    = S \rtensorh T$,
  \item $\alpha_{(\hba, \kdg, \lfe)}$ is the composition of the isomorphism
    \begin{align} \label{eq:rtp-ass-iso}
    (\HadK) \htensor{(\alpha \rt \gamma)}{\phi} L \xrightarrow
    {m_{(\alpha \rt \gamma)} \tl \Id} (\HadK) {_{\rho_{(\alpha \rt \gamma)}}
      \tl} \phi \xrightarrow{(\Id \tr m^{op}_{\delta}) \tl \Id}
    \alpha \tr_{\rho_{\delta}} K {_{\rho_{\gamma}} \tl} \phi
      \end{align}
      with the inverse of the isomorphism
      \begin{align} \label {eq:rtp-ass-iso2}
        H \htensor{\alpha}{(\delta \lt \phi)} (K \htensor{\gamma}{\phi} L) 
        \xrightarrow{\Id \tr m^{op}_{(\delta \lt \phi)}}   \alpha
        \tr_{\rho_{(\delta \lt \phi)}} (K \htensor{\gamma}{\phi} L)
        \xrightarrow{\Id \tr (m_{\gamma} \tl \Id)} \alpha \tr_{\rho_{\delta}} K
        {_{\rho_{\gamma}} \tl } \phi,
      \end{align}
  \item $l_{(\hba)}$ and $r_{(\hba)}$ are given by the compositions
    \begin{align*}
      \frakH \htensor{\frakB}{\beta}  H
      \xrightarrow{m_{\frakB} \tl \Id} \frakH {_{\rho_{\frakB}} \tl} \beta
      \xrightarrow{m^{op}_{\beta}} H, &&
 H \htensor{\alpha}{\frakBo} \frakH 
      \xrightarrow{\Id \tr m_{\frakBo}} \alpha \tr_{\rho_{\frakBo}} \frakH
      \xrightarrow{m_{\alpha}} H.
    \end{align*}
  \end{itemize}
\end{theorem}
\begin{proof}
  Straightforward.
\end{proof}
\begin{remark} \label{remark:cbbimod}
  More explicitly, the isomorphisms \eqref{eq:rtp-ass-iso} and
  \eqref{eq:rtp-ass-iso2} are given by
  \begin{align*}
      [|\alpha\rangle_{1}\gamma] \tr \cbasesb \tl \phi &\to
      \alpha \tr_{\rho_{\delta}} K {_{\rho_{\gamma}} \tl} \phi, &
   \alpha \tr \cbasesb \tl [|\phi\rangle_{2}\delta] &\to \alpha
   \tr_{\rho_{\delta}} K {_{\rho_{\gamma}} \tl} \phi,\\
      |\xi\rangle_{1}\eta \tr \zeta \tl \vartheta &\mapsto \xi \tr \eta\zeta
      \tl \vartheta,   &
   \xi \tr \zeta \tl | \vartheta\rangle_{2}\eta' &\mapsto
   \xi \tr \eta'\zeta \tl \vartheta,
 \end{align*}
 respectively, and $l_{(\hba)}$ and $r_{(\hba)}$ are given by
 \begin{align*}
   \frakB \tr \cbasesb \tl \beta &\to H, \
 b \tr \zeta \tl \xi \mapsto \xi b
   \zeta, &  \alpha \tr \cbasesb \tl \frakBo &\to H, 
  \ \eta \tr \zeta \tl b^{\dag} \mapsto
   \eta b^{\dag} \zeta.
\end{align*}
 \end{remark}

\subsection{The monoidal functor}

\label{subsection:bimod-functor}

We fix the following data:
\begin{itemize}
\item a decomposable $C^{*}$-algebra $B$ with an admissible inverse
  semigroup $\Theta \subseteq \paut(B)$,
\item a covariant representation $(\pi,\upsilon)$ of $(B,\Theta)$ on a
  Hilbert space $\frakK$ (see Section 3), where $\pi$ is faithful.
\end{itemize}
We define $\funct{}$ and $\functs{}$ as in Section 3, and put
  \begin{align*}
    \frakH &:= \funct{B}, & \frakB&:=\functs{\lft{\hmg{B}}}
    \subseteq {\cal L}(\frakH), & \frakBo &:= \functs{\rgt{\hmg{B}}}
    \subseteq {\cal L}(\frakH).
  \end{align*}
\begin{proposition} \label{proposition:cbase}
  $\cbasesb$ is a $C^{*}$-base.
\end{proposition}
\begin{proof}
  Put $\msb:=\hmg{B}$.  By assumption on $\Theta$, the families
  $\lft{\msb}$ and $\rgt{\msb}$ are $\Theta$-supported. By Proposition
  \ref{proposition:functs}, $[\frakB^{*}\frakB] =
  \functs[\lft{\msb}^{*}\lft{\msb}] = \functs{\lft{\msb}} = \frakB$
  and $ [(\frakBo)^{*}\frakBo] = \functs[\rgt{\msb}^{*}\rgt{\msb}] =
  \functs{\rgt{\msb}} = \frakBo$, so $\frakB$ and $\frakBo$ are
  $C^{*}$-algebras. They commute because they are the closed linear
  span of operators of the form $\funct{\lft{b}}$ and
  $\funct{\rgt{c}}$, respectively, where $b \in \hmgs{B}{\theta},c \in
  \hmgs{B}{\theta'}$ for some $\theta,\theta' \in \Theta$ and
  $\funct{(\lft{b})}\funct{(\rgt{c})} =
  \funct{(\lft{b}\rgt{c})}=\funct{(\rgt{c}\lft{b})} =
  \funct{(\lft{b})}\funct{(\rgt{c})}$ for all such $b,c$ by
  Proposition \ref{proposition:funct}.  Finally, $\frakB$ and
  $\frakBo$ act nondegenerately on $\frakH$ because the inclusion
  $Z(B)=\hmgs{B}{\Id} \subseteq B$ is nondegenerate \cite[Proposition
  3.20 (v)]{timmermann:hopf} and $Z(B) \tr_{\pi} \Id \subseteq \frakB
  \cap \frakBo$ acts nondegenerately on $B \opi \frakK=\frakH$.
\end{proof}

We construct a functor from the category of $\Theta$-admissible right
$C^{*}$-$B$-$B$-bimodules to the category of
$C^{*}$-$\cbasesb$-$\cbasesb$-bimodules.  First, we assign to every
$\Theta$-admissible right $C^{*}$-$B$-$B$-bimodule a right
$C^{*}$-$\cbasesb$-$\cbasesb$-bimodule:
\begin{proposition} \label{proposition:modules} Let $E$ be a
  $\Theta$-admissible $C^{*}$-bimodule over $B$.  Then
  \begin{align*}
    \beta(E):= \functs{\rgt{\hmg{E}}} \subseteq {\cal
      L}(\frakH, \funct{E}) \quad \text{and} \quad
    \alpha(E):=\functs{\lft{\hmg{E}}} \subseteq {\cal
      L}(\frakH,\funct{E})
  \end{align*}
  are compatible $C^{*}$-factorizations of $\funct{E}$ with
  respect to $\cbaseosb$ and $\cbasesb$, respectively. The
  representations $\rho_{\beta(E)} \colon \frakB \to {\cal
    L}(\funct{E})$ and $\rho_{\alpha(E)} \colon \frakBo \to {\cal
    L}(\funct{E})$ are given by
  \begin{align} \label{eq:rho-formula}
    \rho_{\beta(E)}(\funct{\lft{b}}) = \funct{\lfte{b}} \quad
    \text{and} \quad
    \rho_{\alpha(E)}(\funct{\rgt{b}}) = \funct{\rgte{b}}
    \end{align}
    for all $b \in \hmgs{B}{\theta}$,  $\theta \in \Theta$.
\end{proposition}
\begin{proof}
  Put $\mse:=\hmg{E}$ and $\msb:=\hmg{B}$.
  First, we show that $\beta(E)$ is a $C^{*}$-factorization of
  $\funct{E}$ with respect to $\cbaseosb$. Using equation
  \eqref{eq:lr-module} and the assumptions on $E$, we find
  \begin{gather*}
    [\beta(E)^{*}\beta(E)] = \functs{[\rgt{\mse}^{*}\rgt{\mse}]} =
    \functs{\rgt{\msb}} = \frakBo, \quad
    [\beta(E)\frakBo] = \functs{[\rgt{\mse}\rgt{\msb}]} = \functs{\rgt{\mse}} = \beta(E), \\
    [\beta(E)\frakH] = \overline{\sum_{\theta \in\Theta}}
    \rgt{\mse_{\theta}}B \opi \upsilon(\theta)\frakK \supseteq
    \overline{\sum_{\theta \in\Theta}} \mse_{\theta} \opi
    \pi(\Img(\theta))\frakK = \overline{\sum_{\theta \in\Theta}}
    \mse_{\theta} \opi \frakK = E \opi \frakK.
  \end{gather*}
  Therefore, $\beta(E)$ is a $C^{*}$-factorization as claimed.  By
  definition of $\beta(E)$ and $\rho_{\beta(E)}$, Lemma
  \ref{lemma:lr-module}, and Proposition \ref{proposition:funct},
  $\rho_{\beta(E)}(\funct{\lft{b}}) \funct{\rgt{\xi}}
    = \funct{(\rgt{\xi})}\funct{(\lft{b})} 
    =\funct{(\lfte{b}})  \funct{(\rgt{\xi})}$
  for all $b \in B$, $\xi \in \mse_{\theta}$, $\theta \in \Theta$.
  This calculation proves the formula for $\rho_{\beta(E)}$ in
  \eqref{eq:rho-formula}.

  Similar calculations show that $\alpha(E)$ is a
  $C^{*}$-factorization of $\funct{E}$ with respect to $\cbasesb$, and
  that $\rho_{\alpha(E)}$ is given by the formula in
  \eqref{eq:rho-formula}.
  
Finally, by equation \eqref{eq:rho-formula}, Proposition
  \ref{proposition:functs}, and Lemma \ref{lemma:lr-module}, 
  \begin{align*}
    [\rho_{\alpha(E)}(\frakBo)\beta(E)] &= \functs{
      [\rgte{\msb}\rgt{\mse}]} = \functs{\rgt{\mse}} = \beta(E), \\
    [\rho_{\beta(E)}(\frakB)\alpha(E)] &=
    \functs{[\lfte{\msb}\lft{\mse}]} = \functs{\lft{\mse}} =
    \alpha(E),
  \end{align*}
  so $\alpha(E)$ and $\beta(E)$ are compatible.
\end{proof}
\begin{remark} \label{remark:alpha} Let $E$ be a right $C^{*}$-module
  over $B$. Then for each $\xi \in E$, there exists an operator
  $\lft{\xi} \in {\cal L}_{B}(B,E)$ such that $\lft{\xi}b=b$ for all
  $b \in B$.  One easily verifies that the map $B \to {\cal
    L}(\funct{B})$ given by $b \mapsto \lft{b}
  \opi \Id$ defines an isomorphism $\Phi_{B} \colon B \to \frakB$ of
  $C^{*}$-algebras. Let $E$ be an  admissible right
  $C^{*}$-$B$-$B$-bimodule and identify $B$ with $\frakB$
  via $\Phi_{B}$. Then the map $E \to {\cal
    L}(\funct{B},\funct{E})$ given by $\xi \mapsto \lft{\xi} \opi \Id$
  defines an isomorphism $\Phi_{E} \colon E \to \alpha(E)$ of
  $C^{*}$-modules over $B \cong \frakB$.
\end{remark}

The next step is to consider morphisms of right $C^{*}$-$B$-$B$-bimodules:
\begin{proposition} \label{proposition:module-morphisms}
  Let $E$ and $F$ be $\Theta$-admissible $C^{*}$-bimodules over
  $B$. Then the assignment $T \mapsto \funct{T}$ defines a bijection of
  ${\cal L}^{B}_{B}(E,F)$ with ${\cal L}\big({_{\beta(E)}\funct{(E)}_{\alpha(E)}},
    {_{\beta(F)}\funct{(F)}_{\alpha(F)}}\big)$.
\end{proposition}
\begin{proof}
  Let $T \in {\cal L}^{B}_{B}(E,F)$. We claim that $\funct{T} \in {\cal
    L}\big({_{\beta(E)}\funct{(E)}_{\alpha(E)}},
  {_{\beta(F)}\funct{(F)}_{\alpha(F)}}\big)$. Indeed, since $T\hmgs{E}{\theta}
  \subseteq \hmgs{F}{\theta}$ for all $\theta \in \Theta$, 
  \begin{align*}
    \funct{(T)}\alpha(E) &\subseteq \functs{[T\lft{\hmg{E}}]} \subseteq
    \functs{\lft{\hmg{F}}} = \alpha(F), \\ \funct{(T)}\beta(E)
    &\subseteq \functs{[T\rgt{\hmg{E}}]} \subseteq \functs{\rgt{\hmg{F}}}=
    \beta(F),
  \end{align*}
  and  similar calculations show that $\funct{(T)}^{*}\alpha(F) \subseteq
  \alpha(E)$ and $\funct{(T)}^{*}\beta(F) \subseteq \beta(E)$.
  
  Since $\pi$ is faithful, the assignment $T \mapsto \funct{T} = T
  \opi \Id$ is injective.

  Finally, let $S \in {\cal
    L}\big({_{\beta(E)}\funct{(E)}_{\alpha(E)}},
  {_{\beta(F)}\funct{(F)}_{\alpha(F)}}\big)$. We show that $S =
  \funct{T}$ for some $T \in {\cal L}^{B}_{B}(E,F)$.  Since
  $S\alpha(E) \subseteq \alpha(F)$ and $S^{*}\alpha(F) \subseteq
  \alpha(E)$, the operator $S$ defines an operator $S_{\alpha(E)} \in
  {\cal L}_{\frakB}(\alpha(E),\alpha(F))$ via $\omega \mapsto
  S\omega$. Put $T:=\Phi_{F}^{-1} S_{\alpha(E)} \Phi_{E} \in {\cal
    L}_{B}(E,F)$, where $\Phi_{F}$ and $\Phi_{E}$ are as in Remark
  \ref{remark:alpha}.  By Lemma \ref{lemma:lr-module},
  \begin{align*}
    \funct{(T)} \funct{(\lft{\xi})} = \funct{(T\lft{\xi})} =
    \funct{(\lft{T\xi})} = S \funct{\lft{\xi}} \quad \text{for all }
    \xi \in \hmgs{E}{\theta}, \theta \in \Theta.
  \end{align*}
  Since $[\functs{(\lft{\hmg{E}})}\frakH]=\funct{E}$, we can conclude
  $\funct{T}=S$.  The assumption on $S$ and equation
  \eqref{eq:trho-commute} imply that for all $b \in \hmgs{B}{\theta}$,
  $\theta \in \Theta$,
    \begin{align*}
      \funct{(T\lfte{b})} = \funct{(T)} \funct{(\lfte{b})} = S
      \rho_{\beta(E)}(\funct{\lft{b}}) =
      \rho_{\beta(F)}(\funct{\lft{b}}) S
      =\funct{(\lftf{b})}\funct{(T)}  = \funct{(\lftf{b}T)}.
    \end{align*}
    Since $\pi$ is injectve and $B$ is $\Theta$-decomposable, we can conclude $T \in {\cal L}^{B}_{B}(E,F)$.
\end{proof}
\begin{corollary}
  There exists a full and faithful functor
  \begin{align*}
    \bbimod&\to \cbbimod
  \end{align*}
  defined on objects by $E \mapsto
  {_{\beta(E)}\funct{(E)}_{\alpha(E)}}$ and on morphisms by $T \mapsto
  \funct{T}$. \qed
\end{corollary}

We show that the functor constructed above is monoidal.
 Let $E$ and $F$ be $\Theta$-admissible right
$C^{*}$-$B$-$B$-bimodules. Then there exists an isomorphism
  \begin{align*}
    \tau_{E,F} \colon \funct{(E)} \htensor{\alpha(E)}{\beta(F)}
    \funct{(F)} & \xrightarrow{\Id \tr m_{\beta(F)}^{op}} \alpha(E)
    \tr_{\rho_{\beta(F)}} \funct{(F)} \xrightarrow{\Phi^{-1}_{E} \tr \Id}
    E \tr F \opi \frakK = \funct{(E \tr F)},
  \end{align*}
where  $\Phi_{E} \colon E \to \alpha(E)$ denotes the isomorphism  $\xi
 \mapsto \lft{\xi} \opi \Id$, see Remark \ref{remark:alpha}.
Explicitly, $\tau_{E,F}$ is given by
\begin{align*}
  \funct{\lft{\xi}} \tr (b \opi \zeta) \tl \funct{\rgt{\eta}}
  \mapsto \xi \tr b\eta \opi 
  \upsilon(\theta')\zeta
\end{align*}
 for all $b \opi\zeta \in   \frakH$, $\xi \in \hmgs{E}{\theta}$,
 $\eta \in \hmgs{F}{\theta'}$, $\theta,\theta' \in \Theta$. Clearly, $\tau_{E,F}$ is
 natural in $E$ and $F$.

 Recall the families $|\hmg{E}\rangle_{1} \subseteq \cL_{\Id}(F,E \tr F)$
 and $|\hmg{F}]_{2} \subseteq \cL^{\Id}(E,E \tr F)$ defined before
 Proposition \ref{proposition:itp}.  Straightforward calculations show
\begin{align}\label{eq:tau-ketbra}
  \begin{aligned}
    \tau_{E,F} |\alpha(E)\rangle_{1} &= \functs{|\hmg{E}\rangle_{1}}
    \subseteq {\cal L}(\funct{F}, \funct{(E \tr F)}), \\ \tau_{E,F}
    |\beta(F)\rangle_{2} &= \functs{|\hmg{F}]_{2}} \subseteq {\cal
      L}(\funct{E}, \funct{(E \tr F)}).
  \end{aligned}
  \end{align}
 \begin{proposition}  \label{proposition:tau-module}
  Let $E$ and $F$ be $\Theta$-admissible right
$C^{*}$-$B$-$B$-bimodules. Then 
\begin{align*}
  \tau_{E,F}(\alpha(E) \rt \alpha(F)) &= \alpha(E \tr F), &
  \tau_{E,F}(\beta(E) \lt \beta(F)) &= \beta(E \tr F).
\end{align*} 
\end{proposition}
\begin{proof}
  Put $G:=E \tr F$.  By equation \eqref{eq:tau-ketbra} and
  Propositions \ref{proposition:functs}, \ref{proposition:itp},
  \begin{gather*}
    \tau_{E,F} (\alpha(E) \rt \alpha(F)) =
    \tau_{E,F}[|\alpha(E)\rangle_{1} \alpha(F)] =
    \functs{[|\hmg{E}\rangle_{1}\lft{\hmg{F}}]} = \functs{\lft{\hmg{G}}} =
    \alpha(G), \\ 
    \tau_{E,F}(\beta(E) \lt \beta(F)) = \tau_{E,F}
    [|\beta(F)\rangle_{2} \beta(E)] = \functs{[|\hmg{F}]_{2}\rgt{\hmg{E}}]} =
    \functs{\rgt{\hmg{G}}} = \beta(G). \qedhere
  \end{gather*}
\end{proof} 
Recall that by definition,
${_{\beta(B)}(\funct{B})_{\alpha(B)}}=\cbaseosb$.
\begin{theorem} \label{theorem:bimodule-functor}
  The full and faithful functor $\bbimod \to \cbbimod$ together with the
  natural transformation $\tau$ defined above and the identity
  $\epsilon:=\Id_{\cbaseosb}$ is a monoidal functor.
\end{theorem}
\begin{proof}
  The main step  is Proposition
  \ref{proposition:tau-module}; the rest  is straightforward but tedious.
\end{proof}

For later use, we note the following analogue of proposition
\ref{proposition:tau-module}:
\begin{proposition} \label{proposition:ket-bra} Let $E$ be a right
  $C^{*}$-module over $B$, let $F$ be $\Theta$-admissible right
  $C^{*}$-$B$-$B$-bimodule, and let $\pi \colon B \to {\cal
    L}_{B}^{B}(F)$ be a nondegenerate representation. Then $E
  \tr_{\pi} F$ is a $\Theta$-admissible right $C^{*}$-$B$-$B$-bimodule
  with respect to the left multiplication given by $b(\xi \tr_{\pi} \eta)=\xi
  \tr_{\pi} b\eta$ for all $b \in B$, $\xi \in E$, $\eta \in F$, and
  \begin{align*}
    \alpha(E \tr_{\pi} F) &= \functs{[|E\rangle_{1}\lft{F}]}, &
    \beta(E \tr_{\pi} F) &= \functs{[|E\rangle_{1}\rgt{F}]}.
  \end{align*}
\end{proposition}
\begin{proof}
  This follows from similar calculations as in the proof of
  Proposition \ref{proposition:tau-module} and from the relation
  $\hmgs{E \tr_{\pi} F}{\theta} = E \tr_{\pi} \hmgs{F}{\theta}$, which
  holds for all $\theta \in \paut(B)$ \cite[Proposition 3.18]{timmermann:hopf}.
\end{proof}
\section{From $C^{*}$-families to $C^{*}$-$\cbasesb$-algebras}

In this section, we extend $\functs{}$ to a functor from the category
of generalized $C^{*}$-algebras, more precisely, the
$\Theta$-admissible $C^{*}$-families used in
\cite{timmermann,timmermann:hopf}, to the category of concrete
$C^{*}$-algebras over $C^{*}$-bases used in \cite{timmer:ckac,timmer:cpmu}. 
Moreover, we show that this functor embeds the fiber product of
$C^{*}$-families into the fiber product of $C^{*}$-algebras.

\subsection{The monoidal category of $\Theta$-admissible $C^{*}$-families}

To define the legs of a pseudo-multiplicative unitary in the form of
Hopf $C^{*}$-bimodules, we introduced in
\cite{timmermann,timmermann:hopf} a monoidal category of generalized
$C^{*}$-algebras: $C^{*}$-families consisting of homogeneous operators
on $C^{*}$-bimodules, morphisms of such $C^{*}$-families, and fiber
products of $C^{*}$-families and morphisms.  In the following
paragraphs, we recall these concepts and introduce the class of
normal morphisms. As before, let $B$ be a fixed $C^{*}$-algebra.

\begin{definition*}[{\cite{timmermann:hopf}}]
  Let $E$ be a right $C^{*}$-$B$-$B$-bimodule.  A {\em $C^{*}$-family}
  on $E$ is a family of closed subspaces $\msc \subseteq \cL(E)$
  satisfying $[\msc^{*}\msc] = \msc$.  We call such a $C^{*}$-family
  $\msc$ {\em nondegenerate} if $[\msc E]=E$, and define its {\em
    multiplier $C^{*}$-family} $\msm(\msc) \subseteq \cL(E)$ by
  \begin{align*}
    \msm(\msc)^{\rho}_{\sigma} := \{ T \in \cL^{\rho}_{\sigma}(E) \mid
    [T\msc], [\msc T] \subseteq \msc\} \quad \text{for all } \rho \in
    \paut(B), \sigma \in \paut(B).
  \end{align*}
 \end{definition*}
%  Let $E$ be a right $C^{*}$-$B$-$B$-bimodule and let $\msc$ be a
%  $C^{*}$-family on $E$.  Then $\msc^{\Id} \subseteq \cL^{\Id}(E)$ and
%  $\msc_{\Id} \subseteq \cL_{\Id}(E)$ are $C^{*}$-families again.
%  Indeed, the assumption on $\msc$ implies
%  $[(\msc^{\Id})^{*}\msc^{\Id}] \subseteq \msc^{\Id}$ and
%  $[(\msc_{\Id})^{*}\msc_{\Id}] \subseteq \msc_{\Id}$, and the reverse
%  inclusion follows from \cite[Remarks 3.9 (iii)]{timmermann:hopf}.

 Let $E$ be a right $C^{*}$-$B$-$B$-bimodule. Then the family
$\mso(E)\subseteq \cL(E)$ given by
\begin{align*}
  \mso(E)^{\rho}_{\sigma} :=
  [\lfte{\hmgs{B}{\rho}}\rgte{\hmgs{B}{\sigma^{*}}}] \quad \text{for
    all } \rho,\sigma \in \paut(B)
\end{align*}
is a $C^{*}$-family \cite[Proposition 3.21]{timmermann:hopf}.

% \begin{proposition}
%   The family of closed subspaces $\cK(B) \subseteq \cL(E)$ given by
% \begin{align*}
%   \cK(B)^{\rho}_{\sigma}:= \{ \theta \circ \lft{\hmgs{B}{\rho'}}
%   \rgt{\hmgs{B}{\sigma'}} \mid \rho' ,\sigma',\theta \in \paut(B),
%   \theta\rho' \leq \rho, \, \theta\sigma'{}^{*} \leq \sigma, \,
%   \Img(\rho') \cap \Img(\sigma') \subseteq \theta \}
% \end{align*}  
% for all $\rho,\sigma \in \paut(B)$ is a $C^{*}$-family.
% \end{proposition}
% \begin{proof}
  
% \end{proof}

\begin{definition}
  Let $\Theta\subseteq \paut(B)$ be an admissible inverse semigroup. A
  $C^{*}$-family $\msc$ on a right $C^{*}$-$B$-$B$-bimodule $E$ is
  {\em $\Theta$-admissible} if $E$ is $\Theta$-admissible, $\msc$ is
  nondegenerate and $\Theta$-supported, $[\msc \mso(E)] \subseteq
  \msc$, and $[\msc^{\rho}_{\sigma}
  \mso^{\rho^{*}\rho}_{\sigma^{*}\sigma}(E)] = \msc^{\rho}_{\sigma}$
  for all $\rho,\sigma \in \Theta$.
\end{definition}
\begin{remark}
  Let $\Theta\subseteq \paut(B)$ be an admissible inverse semigroup
  and $E$ a $\Theta$-admissible right $C^{*}$-$B$-$B$-bimodule. Then
   the $C^{*}$-family $\mso(E)$ is
  $\Theta$-admissible, as one can easily check.
\end{remark}
The following $C^{*}$-family will turn out to be the unit for the
fiber product:
\begin{proposition}
  Let $\Theta \subseteq \paut(B)$ be an admissible inverse
  semigroup. Put $\msb:=\hmg{B}$.
  \begin{enumerate}
  \item Let $\rho,\sigma,\theta \in \paut(B)$, $a \in
    \msb_{\rho}$, $c \in \msb_{\sigma^{*}}$, $d \in
    \msb_{\theta^{*}\theta}$. Then the map
    \begin{align*}
      k^{\theta,d}_{a,c} \colon B \to B, \ b \mapsto \theta(abcd),
    \end{align*}
    is a $(\theta\rho,\theta\sigma)$-homogeneous operator.
  \item The family $\cK(B;\Theta) \subseteq \cL(B)$ given by
    \begin{multline*}
      \cK^{\rho'}_{\sigma'}(B;\Theta) := \clspan \{ k^{\theta,d}_{a,c}
      \mid  a \in \msb_{\rho}, b \in \msb_{\sigma^{*}}, d \in
      \msb_{\theta^{*}\theta},  \rho,\sigma, \theta \in
      \Theta,  \theta\rho \leq \rho', \theta\sigma \leq \sigma'\}
    \end{multline*}
    is a $\Theta$-admissible $C^{*}$-family.
  \end{enumerate}
  \begin{proof}
    For all $\theta \in \paut(B)$ and $d \in
    \msb_{\theta^{*}\theta}$,  define
    $s^{\theta,d} \colon B \to B$ by $b \mapsto \theta(bd)$. 

\smallskip

    i) Since $k^{\theta,d}_{a,c}= s^{\theta,d}\circ \lft{a}\rgt{c}$,
    it suffices to prove $s^{\theta,d} \in
    \cL^{\theta}_{\theta}(B)$.  Evidently, $\Img s^{\theta,d}
    \subseteq \Img(\theta)B$ and
    $s^{\theta,d}(b'b)=\theta(b')s^{\theta,d}(b)$ for all $b' \in
    \Dom(\theta)$, $b \in B$. Moreover, for all $b,b' \in B$,
    \begin{align*}
      \langle b|s^{\theta,d}b'\rangle = b^{*}\theta(b'd) =
      \theta(\theta^{*}(b^{*}\theta(d))b') =
      \theta(\theta^{*}(b\theta(d^{*}))^{*}b') = \theta(\langle
      s^{\theta^{*},\theta(d^{*})}b|b'\rangle);
    \end{align*}
    here, we used the fact that $d$ is central. Note that
    $\theta(d^{*})\in \msb_{\theta\theta^{*}}$ by \cite[Proposition
    3.20]{timmermann:hopf}.
    
    \smallskip
    
    ii) Let $a,c,d$ as in i), where  $\rho,\sigma,\theta \in \Theta$. Write
    $d=d_{a}d'd_{c}$ with $d_{a},d',d_{c} \in \msb_{\theta^{*}\theta}$
    and put $a':=\theta^{*}(a\theta(d_{a}))$ and
    $c':=\theta^{*}(\theta(d_{c})c)$. Then $a' \in
    \msb_{\theta^{*}\rho\theta}$, $c' \in
    \msb_{\theta^{*}\sigma^{*}\theta}$, and
    \begin{align*}
      a\theta(bd)c =a\theta(d_{a})\theta(bd')\theta(d_{c})c
      =\theta(a'bd'c') \quad \text{for all } b \in B.
    \end{align*}
    Therefore, $\lft{a}\rgt{c} \circ s^{\theta,d} = s^{\theta,d'}
    \circ \lft{a'}\rgt{c'}$. Combining this observation with the
    results from (the proof of) i), we can conclude that
    $\cK(B;\Theta)$ is a $C^{*}$-family. By definition,
    $\cK(B;\Theta)$ is $\Theta$-supported. Using the facts that $B$ is
    $\Theta$-decomposable, $\Id_{B} \in \Theta$, and that $\mso(B)$ is
    $\Theta$-supported, one easily verifies that the $C^{*}$-family
    $\cK(B;\Theta)$ is nondegenerate and satisfies $[\cK(B;\Theta)
    \mso(B)] = \cK(B;\Theta)$.

    Finally, if $a,b,c,d$ and $d_{a},d',d_{c}$ are as above, then
    $k^{\theta,d}_{a,c} = k^{\theta,d'}_{d_{a}a,cd_{c}} =
    s^{\theta,d'} \circ \lft{d_{a}a}\rgt{cd_{c}}$ and
    $\lft{d_{a}a}\rgt{cd_{c}} \in \mso^{\rho'}_{\sigma'}(B)$, where
    $\rho'=\theta^{*}\theta\rho$ and $\sigma'=\theta^{*}\theta\sigma$.
    By \cite[Remark 3.9 (ii)]{timmermann:hopf}, $\mso^{\rho'}_{\sigma'}(B) =
    [\mso^{\rho'}_{\sigma'}(B)
    \mso^{\rho'{}^{*}\rho'}_{\sigma'{}^{*}\sigma'}(B)]$, and hence
    \begin{align*}
      k^{\theta,d}_{a,c} \in
      [\cK^{\theta\rho}_{\theta\sigma}(B;\Theta)
      \mso^{\rho'{}^{*}\rho'}_{\sigma'{}^{*}\sigma'}(B)]=
      [\cK^{\theta\rho}_{\theta\sigma}(B;\Theta)
      \mso^{\rho^{*}\theta^{*}\theta\rho}_{\sigma^{*}\theta^{*}\theta\sigma}(B)].
    \end{align*}
    Consequently, $\cK^{\rho}_{\sigma}(B;\Theta) \subseteq
    [\cK^{\rho}_{\sigma}(B;\Theta)
    \mso^{\rho{}^{*}\rho}_{\sigma{}^{*}\sigma}(B)]$ for all
    $\rho,\sigma \in \Theta$.
  \end{proof} 
\end{proposition}

The fiber product of $C^{*}$-families is defined as follows.  Let $E$
and $F$ be right $C^{*}$-$B$-$B$-bimodules and let $\msc \subseteq
\cL(E)$ and $\msd \subseteq \cL(F)$ be families of closed subspaces.
We call two partial automorphisms $\rho,\sigma \in \paut(B)$ {\em
  compatible} and write $\rho \perp \sigma$ if $\rho\sigma^{*} \leq
\Id$ and $\rho^{*}\sigma \leq \Id$.  By \cite[Proposition
5.3]{timmermann:hopf}, there exists for all $\rho,\rho',\sigma,\sigma'
\in \paut(B)$ satisfying $\sigma' \perp \rho'$ a map
\begin{align*}
  \cL^{\rho}_{\sigma'}(E) \times \cL^{\rho'}_{\sigma} (F) \to
  \cL^{\rho}_{\sigma} (E \tr F), \quad (S,T) \mapsto S
  \tr T,
\end{align*}
where $(S \tr T) (\xi \tr \eta) = S\xi \tr T\eta$ for all $S \in
\cL^{\rho}_{\sigma'}(E)$, $T \in \cL^{\rho'}_{\sigma}(F)$ and $\xi \in
E$, $\eta \in F$. 

We define a family of closed subspaces $\msc \tr \msd \subseteq \cL(E
\tr F)$ by
  \begin{align*}
    (\msc \tr \msd)^{\rho}_{\sigma} := \clspan \{ c \tr d \mid
    c \in \msc^{\rho}_{\sigma'}, d \in \msd^{\rho'}_{\sigma}, \,
    \sigma',\rho' \in \paut(B) \text{ compatible}\}
  \end{align*}
  for all $\rho,\sigma \in \paut(B)$.  If $\msc$ and $\msd$ are
  (nondegenerate) $C^{*}$-families, then so is $\msc \tr \msd$, as one
  can easily check. Moreover, in that case, $\msm(\msc) \tr \msm(\msd)
  \subseteq \msm(\msc \tr \msd)$. 
  \begin{proposition}
    Let $\Theta \subseteq \paut(B)$ be an admissible inverse
    semigroup.  If $\msc$ and $\msd$ are $\Theta$-admissible
    $C^{*}$-families, then so is $\msc \tr \msd$. 
  \end{proposition}
  \begin{proof}
    Let $\msc$ and $\msd$ be $\Theta$-admissible $C^{*}$-families on
    right $C^{*}$-$B$-$B$-bimodules $E$ and $F$, respectively, and put
    $\msb:=\hmg{B}$. Then $E
    \tr F$ is $\Theta$-admissible by Proposition
    \ref{proposition:itp}, $\msc \tr \msd$ is nondegenerate and
    $\Theta$-supported by construction and assumption of $\msc$ and
    $\msd$, and $[(\msc \tr \msd)\mso(E \tr F)] = [\msc
    \lfte{\msb}] \tr [\msd \rgtf{\msb}] \subseteq \msc \tr
    \msd$. Finally, by assumption on $\msc$ and $\msd$, 
    \begin{align*}
      [(\msc \tr
      \msd)^{\rho}_{\sigma}\mso^{\rho^{*}\rho}_{\sigma^{*}\sigma}(E
      \tr F)] = [(\msc \tr \msd)^{\rho}_{\sigma}
      (\lfte{\msb_{\rho^{*}\rho}} \tr \rgtf{\msb_{\sigma^{*}\sigma}})]
      = (\msc \tr \msd)^{\rho}_{\sigma}
    \end{align*}
    for each $\rho,\sigma \in \Theta$.
  \end{proof}
In \cite{timmermann,timmermann:hopf}, we introduced a rather unwieldy
notion of morphisms between $C^{*}$-families. For our purposes, it
suffices to consider the following special class of morphisms:
\begin{definition}
  Let $\msc$ and $\msd$ be $C^{*}$-families on right
  $C^{*}$-$B$-$B$-bimodules $E$ and $F$, respectively. A  {\em normal
    morphism} from $\msc$ to $\msd$ is a family of maps
  $\phi^{\rho}_{\sigma} \colon \msc^{\rho}_{\sigma} \to
  \msd^{\rho}_{\sigma}$, given for all $\rho,\sigma \in \paut(B)$,
  such that $[I_{\phi}E]=F$, where
    \begin{multline} \label{eq:morphism-intertwiners}
      I_{\phi}:=\big\{ T \in {\cal L}^{B}_{B}(E,F) \,\big|\, Tc =
      \phi^{\rho}_{\sigma}(c)T \text{ and } c T^{*} =
      T^{*}\phi^{\rho}_{\sigma}(c)\\ \text{ for all } c \in
      \msc^{\rho}_{\sigma} \text{ and }\rho,\sigma \in \paut(B) \big\}.
    \end{multline}
    A normal morphism $\phi$ is {\em nondegenerate} if
    $[\phi(\msc)\msd] = \msd$, where $\phi(\msc)=
    ([\phi^{\rho}_{\sigma}(\msc^{\rho}_{\sigma})])_{\rho,\sigma}$.
\end{definition}
Evidently, the composition of normal morphisms is a normal morphism
again. Moreover,  normal morphisms preserve all structure maps:
\begin{proposition} \label{proposition:family-morphism} Let $\msc$ and
  $\msd$ be $C^{*}$-families on right $C^{*}$-$B$-$B$-bimodules $E$
  and $F$, respectively, and let
  $\phi=(\phi^{\rho}_{\sigma})_{\rho,\sigma}$ be a normal morphism
  from $\msc$ to $\msd$. 
  \begin{enumerate}
  \item For all $\rho,\rho',\sigma,\sigma' \in \paut(B)$, $c,c'' \in
    \msc^{\rho}_{\sigma}$, $c' \in \msc^{\rho'}_{\sigma'}$ and
    $\lambda,\mu \in \complex$
    \begin{gather*}
      \phi^{\rho}_{\sigma}(\lambda c + \mu c'') = \lambda
      \phi^{\rho}_{\sigma}(c) + \mu \phi^{\rho}_{\sigma}(c''), \quad
      \phi^{\rho}_{\sigma}(c)\phi^{\rho'}_{\sigma'}(c') =
      \phi^{\rho\rho'}_{\sigma\sigma'}(cc'), \quad
      \big(\phi^{\rho}_{\sigma}(c)\big)^{*} =
      \phi^{\rho^{*}}_{\sigma^{*}}(c^{*}),
      \\
      \phi^{\rho}_{\sigma}(c) = \phi^{\rho'}_{\sigma'}(c) \quad
      \text{if } \rho \leq \rho', \sigma \leq \sigma'.
    \end{gather*}
    In particular, $\phi^{\Id}_{\Id}\colon \msc^{\Id}_{\Id} \to
    \msd^{\Id}_{\Id}$ is a $*$-homomorphism of $C^{*}$-algebras, and
    $\phi(\msc)$ is a $C^{*}$-family.
  \item If $\msc \subseteq \cL(E)$ is nondegenerate, then also
    $\phi(\msc) \subseteq \cL(F)$ is nondegenerate.
  \item If $[\msc\mso(E)] \subseteq \msc$, then for all
    $\rho,\sigma,\theta \in \paut(B)$, $c \in \msc^{\rho}_{\sigma}$,
    $b \in \hmgs{B}{\theta}$,
  \begin{gather*}
    \phi^{\rho}_{\sigma\theta^{*}}(c\rgte{b}) =
    \phi^{\rho}_{\sigma}(c) \rgtf{b}, \qquad
    \phi^{\rho}_{\sigma\theta}(c\lfte{b}) =
    \phi^{\rho}_{\sigma}(c) \lftf{b};
  \end{gather*}
  in particular, $[\phi(\msc)\mso(F)]\subseteq \phi(\msc)$.
\item If $\msc$ and $F$ are $\Theta$-admissible, then $\phi(\msc)$ is 
  $\Theta$-admissible.
  \end{enumerate}
\end{proposition}
\begin{proof}
  i) The equations follow from the facts that $[I_{\phi}E]=F$ and
  that for all $T \in I_{\phi}$,
  \begin{gather*}
    \phi^{\rho}_{\sigma}(\lambda c+\mu c'') T=    T(\lambda c+\mu c'')
    = \lambda T c + \mu T c''  = (\lambda\phi^{\rho}_{\sigma}(c) +
    \mu\phi^{\rho}_{\sigma}(c''))T, \\
    \phi^{\rho\rho'}_{\sigma\sigma'}(cc') T = Tcc' =
    \phi^{\rho}_{\sigma}(c)Tc' =
    \phi^{\rho}_{\sigma}(c)\phi^{\rho'}_{\sigma'} (c') T, \\
    (\phi^{\rho}_{\sigma}(c))^{*}T =
    (T^{*}\phi^{\rho}_{\sigma}(c))^{*} = (cT^{*})^{*} = Tc^{*} =
    \phi^{\rho^{*}}_{\sigma^{*}}(c^{*}) T. 
  \end{gather*}
  
  ii) Evident from the relation $[\phi(\msc)F] = [\phi(\msc)I_{\phi}E]
  = [I_{\phi}\msc E] = [I_{\phi}E]=F$.
  
  \smallskip
  
  iii) The equations follow from the facts that $[I_{\phi}E]=F$ and
  that for all $T \in I_{\phi}$,
  \begin{gather*}
    \phi^{\rho}_{\sigma\theta^{*}}(c\rgte{b})T=Tc\rgte{b} =
    \phi^{\rho}_{\sigma}(c)T\rgte{b}=
    \phi^{\rho}_{\sigma}(c)\rgtf{b}T, \\
    \phi^{\rho}_{\sigma\theta}(c\lfte{b})T=Tc\lfte{b} =
    \phi^{\rho}_{\sigma}(c)T\lfte{b}=
    \phi^{\rho}_{\sigma}(c)\lftf{b}T.
  \end{gather*}

  iv) This follows easily from i)-iii).
\end{proof}

 Normal morphisms are morphisms in the sense of
\cite{timmermann,timmermann:hopf}:
\begin{proposition} \label{proposition:morphism-extend} Let $\msc$ and
  $\msd$ be $C^{*}$-families on right $C^{*}$-$B$-$B$-bimodules $E$
  and $F$, respectively, and let $\phi$ be a normal morphism from
  $\msc$ to $\msd$.  Then for each right
  $C^{*}$-$\complex$-$B$-bimodule $X$ and each right
  $C^{*}$-$B$-$\complex$-bimodule $Y$, there exists a unique
  $*$-homomorphism
  \begin{align*}
    \phi^{X}_{Y} \colon \big(\cL(X) \tr \msc \tr
    \cL(Y)\big)^{\Id}_{\Id} \to \big(\cL(X) \tr \msd \tr
    \cL(Y)\big)^{\Id}_{\Id}
  \end{align*}
  such that for all $U \in \cL^{\Id}_{\sigma}(X)$, $c \in
  \msc^{\rho'}_{\sigma'}$,  $V \in \cL^{\rho}_{\Id}(Y)$, 
  $\sigma,\sigma',\rho,\rho' \in \paut(B)$, where $\sigma \perp \rho'$,
  $\sigma' \perp \rho$,
  \begin{align} \label{eq:morphism}
    \phi^{X}_{Y}(U \tr c \tr V) = U \tr \phi^{\rho'}_{\sigma'}(c) \tr V.
  \end{align}
\end{proposition}
\begin{proof}
  Let $X,Y$ as above, put $I:=\Id_{X} \tr I_{\phi} \tr \Id_{Y}$, and
  let $R \in \big(\cL(X) \tr \msc \tr \cL(Y)\big)^{\Id}_{\Id}$.  Since
  every element of $I_{\phi}^{*}I_{\phi}$ commutes with every element
  of $\msc$, every element of $I^{*}I$ commutes with $R$, and $\langle
  S\upsilon | T R \omega\rangle = \langle \upsilon | S^{*}T R
  \omega\rangle = \langle \upsilon | R S^{*}T \omega\rangle = \langle
  S R^{*}\upsilon | T \omega\rangle$ for all $S,T \in I$ and
  $\upsilon,\omega \in X \tr E \tr Y$. By assumption on $\phi$,
  elements of the form $S\upsilon$ and $T\omega$ as above are linearly
  dense in $X \tr F \tr Y$.  Therefore, the maps
    \begin{align*}
      \phi^{X}_{Y}(R)^{*} \colon X \tr F \tr Y \to X \tr F \tr Y, \quad
      S\upsilon \mapsto SR^{*}\upsilon, \\
      \phi^{X}_{Y}(R) \colon X \tr F \tr Y \to X \tr F \tr Y, \quad
      T\omega \mapsto TR\omega,
    \end{align*}
    where $S,T \in I$ and $\upsilon,\omega \in X \tr E \tr Y$, are
    well-defined and form an adjoint pair of operators.  By
    definition, equation \eqref{eq:morphism} holds; in particular,
    $\phi^{X}_{Y}((\cL(X) \tr \msc \tr \cL(Y))^{\Id}_{\Id}) \subseteq
    (\cL(X) \tr \msd \tr \cL(Y))^{\Id}_{\Id}$.
\end{proof}
\begin{proposition} \label{proposition:morphism-multiplier}
  Let $\msc$ and $\msd$ be nondegenerate $C^{*}$-families on right
  $C^{*}$-$B$-$B$-bimodules $E$ and $F$, respectively, and let $\phi$ be a
  normal nondegenerate morphism from $\msc$ to $\msm(\msd)$. Then
  $\phi$ extends uniquely to a normal morphism $\tilde \phi$ from
  $\msm(\msc)$ to $\msm(\msd)$. 
\end{proposition}
\begin{proof}
  The proof is similar to the proof above.  Let $R \in
  \msm(\msc)^{\rho}_{\sigma}$, where $\rho,\sigma \in \paut(B)$.
  Since every element of $I_{\phi}^{*}I_{\phi}$ commutes with every
  element of $\msc$ and $\msc$ acts nondegenerately on $E$, every
  element of $I_{\phi}^{*}I_{\phi}$ commutes with $R$. Therefore,
  $\langle S \upsilon| T R\omega\rangle = \sigma(\langle S R^{*}\upsilon|
  T\omega\rangle)$ for all $S,T \in I_{\phi}$ and $\upsilon,\omega \in
  E$. Since $[I_{\phi}E]=E$, the maps
  \begin{align*}
    \tilde \phi^{\rho}_{\sigma}(R)^{*} \colon F \to F, \ S\upsilon
    \mapsto S R^{*}\upsilon, \quad \text{and} \quad \tilde
    \phi^{\rho}_{\sigma}(R) \colon F \to F, \ T \mapsto TR\omega,
  \end{align*}
  where $S,T \in I_{\phi}$ and $\upsilon,\omega \in E$, are
  well-defined. One easily checks that $ \tilde
  \phi^{\rho}_{\sigma}(R)^{*} \in \cL^{\rho^{*}}_{\sigma^{*}}(F)$ and
  $\phi^{\rho}_{\sigma}(R) \in \cL^{\rho}_{\sigma}(F)$. 
  
  Letting $\rho,\sigma$, $R$ vary, we obtain a family
  of maps $(\tilde \phi^{\rho}_{\sigma})_{\rho,\sigma}$, where each $\tilde
  \phi^{\rho}_{\sigma}$ extends $\phi^{\rho}_{\sigma}$. Since $\phi$
  is nondegenerate, $[\tilde\phi(\msm(\msc))\msd] = [\tilde
  \phi(\msm(\msc))\phi(\msc)\msd] = [\tilde\phi(\msm(\msc)\msc)\msd] =
  [\phi(\msc)\msd]=\msd$, so that $\tilde \phi(\msm(\msc))\subseteq
  \msm(\msd)$. Clearly, $I_{\phi} \subseteq I_{\tilde \phi}$, whence
  $\tilde \phi$ is a normal morphism.
\end{proof}

The fiber product is functorial in the following sense:
  \begin{proposition} \label{proposition:morphism-tensor}
    Let $\phi \colon \msa \to \msb$ and $\psi \colon \msc \to \msd$ be
    normal morphisms of $C^{*}$-families on right
    $C^{*}$-$B$-$B$-bimodules. Then there exists a unique normal morphism
    $\phi \tr \psi \colon \msa \tr \msc \to \msb
    \tr \msd$ such that 
    \begin{align} \label{eq:morphism-tensor}
      (\phi \tr \psi)^{\rho}_{\sigma}(a \tr c) =
      \phi^{\rho}_{\sigma'}(a) \tr \phi^{\rho'}_{\sigma}(c) 
    \end{align}
    for all $a \in \msa^{\rho}_{\sigma'}$, $c \in
    \msc^{\rho'}_{\sigma}$,  $\rho,\rho',\sigma,\sigma' \in
    \paut(B)$, where $\sigma' \perp \rho'$.
  \end{proposition}
  \begin{proof}
    We follow the same scheme as in the proofs of Propositions
    \ref{proposition:morphism-extend} and
    \ref{proposition:morphism-multiplier}.  Denote by $E$ and $F$ the
    underlying $C^{*}$-bimodules of $\msa \tr \msc$ and $\msb
    \tr \msd$, respectively, and put $I:=I_{\phi} \tr
    I_{\psi} \subseteq {\cal L}^{B}_{B}(E,F)$.

    Let $\rho,\sigma \in \paut(B)$ and $R \in (\msa \tr
    \msc)^{\rho}_{\sigma}$. Since every element of $I^{*}I$ commutes
    with $R$, we have $\langle S\upsilon| T R\omega\rangle =
    \sigma(\langle SR^{*}\upsilon|T\omega\rangle)$ for all $S,T \in I$
    and $\upsilon,\omega \in E$. Since  $[IE]=F$ by assumption on
    $\phi$ and $\psi$, the maps
    \begin{align*}
       (\phi \tr \psi)^{\rho}_{\sigma}(R)^{*} \colon F \to F, \
      S\upsilon \mapsto SR^{*}\upsilon, \quad \text{and} \quad
      (\phi \tr \psi)^{\rho}_{\sigma}(R) \colon F \to F, \
      T\omega \mapsto TR\omega,
    \end{align*}
    are well-defined, and one easily checks that $(\phi \tr
    \psi)^{\rho}_{\sigma}(R)^{*} \in \cL^{\rho^{*}}_{\sigma^{*}}(F)$,
    $(\phi \tr \psi)^{\rho}_{\sigma}(R) \in
    \cL^{\rho}_{\sigma}(F)$.

    Letting $\rho,\sigma$, $R$ vary, we obtain a family of maps
    $((\phi \tr \psi)^{\rho}_{\sigma})_{\rho,\sigma}$.  One
    readily verifies that equation \eqref{eq:morphism-tensor} holds;
    in particular, the image of $(\phi \tr \psi)^{\rho}_{\sigma}$
    is contained in $\msb \tr \msd$. By definition,
    $I \subseteq I_{\phi \tr \psi}$, whence $\phi \tr
    \psi$ is a normal morphism.
  \end{proof}
  Let $\phi \colon \msa \to \msm(\msb)$ and $\psi \colon \msc \to
  \msm(\msd)$ be nondegenerate normal morphism of nondegenerate
  $C^{*}$-families on right $C^{*}$-$B$-$B$-bimodules. Then we obtain a
  nondegenerate normal morphism
  \begin{align*}
    \phi \tr \psi \colon \msa
  \tr \msc \to \msm(\msb) \tr \msm(\msd)
  \hookrightarrow \msm(\msb \tr \msd),
  \end{align*} 
  which we denote by $\phi \tr \psi$ again.

  Let $\Theta \subseteq \paut(B)$ be an admissible inverse semigroup.
  We denote by $\bfamily$ the category whose objects are all
  $\Theta$-admissible $C^{*}$-families on right
  $C^{*}$-$B$-$B$-bimodules and whose morphisms from a family $\msc$
  to a family $\msd$ are all nondegenerate normal morphisms from
  $\msc$ to $\msm(\msd)$.

  \begin{theorem}
    The category $\bfamily$ carries a structure of a monoidal
    category, where $\cK(B;\Theta)$ is the unit and for all  $\Theta$-admissible
    $C^{*}$-families  $\msa,\msc,\msd$ on $C^{*}$-$B$-$B$-bimodules $E,F,G$
    and all  morphisms $\phi,\psi$,
    \begin{itemize}
    \item  $\msc \odot \msd= \msc \tr \msd$ and $\phi \odot \psi= \phi \tr \psi$,
    \item the isomorphisms
    \begin{gather*}
      \alpha_{\msa,\msc,\msd} \colon (\msa \tr \msc)
      \tr \msd \to \msa \tr (\msc \tr \msd),
      \\
      l_{\msc} \colon \cK(B;\Theta) \tr \msc \to \msc, \qquad
      r_{\msc}  \colon \msc \tr \cK(B;\Theta) \to \msc
    \end{gather*}
    are given by conjugation by the isomorphisms
    \begin{align*}
      \alpha_{E,F,G} \colon (E \tr F) \tr G &\to E
      \tr (G \tr G), & l_{E}\colon B \tr E &\to E,
      & r_{E} \colon E \tr B &\to E.
    \end{align*}
    \end{itemize}
  \end{theorem}
  \begin{proof}
    Almost all details that have to be checked are
    straightforward; we only prove that for each $\Theta$-admissible
    $C^{*}$-family $\msc$ on a right $C^{*}$-$B$-$B$-bimodule $E$, the
    morphism $l_{\msc}$ is an isomorphism.  Put $\msb:=\hmg{B}$.

    Let $\rho,\sigma \in \Theta$ and $x \in \msc^{\rho}_{\sigma}$.
    Since $\msc^{\rho}_{\sigma} =
    [\msc^{\rho}_{\sigma}\mso^{\rho^{*}\rho}_{\sigma^{*}\sigma}(E)]$,
    we can write $x=x'\lfte{a}$ with $x' \in \msc^{\rho}_{\sigma}$ and
    $a \in \msb_{\rho^{*}\rho}$. Then the operator $s^{\rho,a} \colon
    B \to B$ given by $b \mapsto \rho(ba)$ belongs to
    $\cK^{\rho}_{\rho}(B;\Theta)$, and
    \begin{align*}
      x l_{E}(b \tr \xi) = xb\xi = x'ab\xi = \rho(ab)x'\xi =
      l_{E}(s^{\rho,a}b \tr x'\xi) \quad \text{for all } b\in B, \xi
      \in E,
    \end{align*}
    showing that $x = \Ad_{l_{E}}(s^{\rho,a} \tr x') \in
    \Ad_{l_{E}}((\cK(B;\Theta) \tr \msc)^{\rho}_{\sigma})$. Since
    $\msc$ is $\Theta$-supported, we can conclude $\msc \subseteq
    \Ad_{L_{E}}(\cK(B;\Theta) \tr \msc)$.

    Conversely, let $x$ be as above and let $\rho',\sigma',\theta \in
    \Theta$, $a \in \msb_{\rho'}$, $c \in \msb_{\sigma'{}^{*}}$, $d
    \in \msb_{\theta^{*}\theta}$ such that $\rho \perp \theta\sigma'$. Write
    $d=d_{a}d_{c}$ with $d_{a},d_{c} \in
    \msb_{\theta^{*}\theta}$. Then $bcd_{c}=\sigma'(cd_{c}b)$ for
    all $b \in B$ by \cite[Proposition
    3.20]{timmermann:hopf}, and hence
    \begin{align*}
      l_{E}(k^{\theta,d}_{a,c}b \tr x\xi) = \theta(abcd)x\xi &=
      \theta(ad_{a}) \theta(\sigma'(cd_{c}b))x\xi \\ &=
      \theta(ad_{a})xcd_{c}b\xi = \lfte{\theta(ad_{a})}x\lfte{cd_{c}}l_{E}(b \tr\xi)
    \end{align*}
    for all $b\in B$, $\xi \in E$. Since $[\msc\mso(E)] \subseteq
    \msc$, $\theta(ad_{a}) \in \msb_{\theta\rho'\theta^{*}}$, $cd_{c}
    \in \msb_{\sigma'{}^{*}\theta^{*}\theta}$, and
    $\rho\sigma'{}^{*}\theta^{*}\theta \leq \theta$, we have
    $\lfte{\theta(ad_{a})}x\lfte{cd_{c}} \in
    \msc^{\theta\rho'\theta^{*}\rho\sigma'{}^{*}\theta^{*}\theta}_{\sigma}
    \subseteq \msc^{\theta\rho'}_{\sigma}$. Consequently, $
    \Ad_{L_{E}}(\cK(B;\Theta) \tr \msc) \subseteq \msc$.
  \end{proof}

  \subsection{The  category of concrete
    $C^{*}$-$\cbasesb$-$\cbaseosb$-algebras}

\label{subsection:cbalg}

  In \cite{timmer:cpmu}, we introduced a fiber product for certain
  $C^{*}$-algebras represented on $C^{*}$-bimodules, and used this
  construction to define reduced Hopf $C^{*}$-bimodules. Let us review
  the pertaining definitions. Throughout this subsection, let
  $\cbasesb$ be a $C^{*}$-base.

  \begin{definition*}[{\cite{timmer:ckac,timmer:cpmu}}]
    A {\em (nondegenerate) concrete $C^{*}$-$\cbasesb$-algebra}
    $(H,A,\alpha)$, briefly written $(H_{\alpha},A)$, consists of a
    Hilbert space $H$, a (nondegenerate) $C^{*}$-algebra $A \subseteq
    {\cal L}(H)$, and a $C^{*}$-factorization $\alpha \in
    \cfact(H;\cbasesb)$ such that $\rho_{\alpha}(\frakBo)A \subseteq
    A$.  If $(H_{\alpha},A)$ is a nondegenerate concrete
    $C^{*}$-$\cbasesb$-algebra, then $A' \subseteq
    \rho_{\alpha}(\frakBo)'$.

    Let $(H_{\alpha},A)$ and $(K_{\beta},B)$ be concrete
    $C^{*}$-$\cbasesb$-algebras. A {\em morphism} from $(H_{\alpha},A)$
    to $(K_{\beta},B)$ is a $*$-homomorphism $\pi \colon A \to B$ such
    that $\beta = [I_{\pi}\alpha]$, where
    \begin{align*}
      I_{\pi} := \big\{ T \in {\cal L}(H_{\alpha},K_{\beta}) \,\big|
      \pi(a)T = T a \text{ for all } a \in A\big\}.
    \end{align*}
    Assume that  $(K_{\beta},B)$ is nongegenerate, so that
    $(K_{\beta},M(B))$ is a concrete $C^{*}$-$\cbasesb$-algebra. We
    call a morphism $\phi$ from $(H_{\alpha},A)$ to $(K_{\beta},M(B))$
    {\em nondegenerate} if $[\phi(A)B]=B$.

    The {\em fiber product} of a concrete $C^{*}$-$\cbasesb$-algebra
    $(H_{\alpha},A)$ and a concrete $C^{*}$-$\cbaseosb$-algebra
    $(K_{\delta},B)$ is the $C^{*}$-algebra
    \begin{align*}
      \AfibreB := \big\{ T \in {\cal L}(\HfibreK) \,\big|\, &T\kalpha{1},
      T^{*}\kalpha{1} \subseteq [ \kalpha{1} B]\text{ and }
      T\kdelta{2}, T^{*}\kdelta{2} \subseteq [\kdelta{2}A]\big\}.
    \end{align*}
  \end{definition*}

  Apart from special cases, we do not know whether the fiber
    product $\AfibreB$ of a nondegenerate concrete
  $C^{*}$-$\cbasesb$-algebra $(H_{\alpha},A)$ and a nondegenerate
  concrete $C^{*}$-$\cbaseosb$-algebra $(K_{\delta},B)$ is
  nondegenerate again. If it is nondegenerate,  then $M(A)
  \hfibre{\alpha}{\beta}M(B) \subseteq M(A \hfibre{\alpha}{\beta}
  B)$; see
  \cite[Lemma 2.5]{timmer:ckac}.

 Let $\phi$ be a morphism of nondegenerate concrete
 $C^{*}$-$\cbasesb$-algebras $(H_{\alpha},A)$ and $(L_{\gamma},C)$, and
 let $\psi$ be a morphism of nondegenerate concrete
 $C^{*}$-$\cbaseosb$-algebras $(K_{\beta},B)$ and $(M_{\delta},D)$.  Then
 there exists a unique $*$-homomorphism
 \begin{align*}
  \phi \ast \psi \colon
 \AfibreB \to C \fibre{\gamma}{\frakH}{\delta} D 
 \end{align*}
 such that $(\phi \ast \psi)(T) \cdot (X \rtensorh Y) = (X \rtensorh
 Y) \cdot T$ for all $T \in \AfibreB$, $X \in I_{\phi}$, $Y \in
 I_{\psi}$; see \cite[Proposition 3.13]{timmer:cpmu}.

 Let $(H_{\alpha},A)$, $(K_{\beta},B)$ be nondegenerate concrete
 $C^{*}$-$\cbasesb$-algebras and let $\pi$ be a nondegenerate morphism
 from $(H_{\alpha},A)$ to $(K_{\beta},M(B))$. Then the unique strictly
 continuous extension $\tilde \pi \colon M(A) \to M(B)$ of $\pi$ is a
 morphism from $(H_{\alpha},M(A))$ to $(K_{\beta},M(B))$; see \cite[Lemma
 2.4]{timmer:ckac}.

 \begin{definition*}[{\cite{timmer:ckac,timmer:cpmu}}]
   A {\em (nondegenerate) concrete
     $C^{*}$-$\cbasesb$-$\cbasesb$-algebra} is a pair $(\hba,A)$
   consisting of a $C^{*}$-$\cbasesb$-$\cbasesb$-bimodule $\hba$ and a
   nondegenerate $C^{* }$-algebra $A \subseteq {\cal L}(H)$ such that
   $\rho_{\alpha}(\frakBo)A \subseteq A$ and $\rho_{\beta}(\frakB)A
   \subseteq A$.

   Let $(\hba,A)$ and $(\kdg,B)$ be concrete
   $C^{*}$-$\cbasesb$-$\cbasesb$-algebras.

   A {\em morphism} from $(\hba,A)$ to $(\kdg,B)$ is a
   $*$-homomorphism $\pi \colon A \to B$ such that
   $\delta=[I_{\pi}\beta]$ and $\gamma=[I_{\pi}\alpha]$, where
   \begin{align} \label{eq:morphism-intertwiners-2} I_{\pi} := \big\{
     T \in {\cal L}(\hba,\kdg) \,\big|\, Ta = \pi(a) T \text{ for all
     } a \in A\big\}.
   \end{align}

   The {\em fiber product} of $(\hba,A)$ and $(\kdg,B)$ is the
   $C^{*}$-algebra
   \begin{align*}
     \AfibreB := \big\{ T \in {\cal L}(\HfibreK) \,\big|\,
     &T\kalpha{1}, T^{*}\kalpha{1} \subseteq [ \kalpha{1} B] 
     \text{ and } T\kdelta{2}, T^{*}\kdelta{2} \subseteq
     [\kdelta{2}A]\big\}.
   \end{align*}
\end{definition*}

   Let $(\hba,A)$ and $(\kdg,B)$ be concrete
   $C^{*}$-$\cbasesb$-$\cbasesb$-algebras.
   If $\pi$ is a morphism from $(\hba,A)$ to $(\kdg,B)$, then by
   \cite[Lemma 2.2]{timmer:ckac},
   \begin{align*}
     \pi(a\rho_{\alpha}(b^{\dag}))=\pi(a)\rho_{\gamma}(b^{\dag}) \quad
     \text{and} \quad \pi(a\rho_{\beta}(b))=\pi(a)\rho_{\delta}(b)
     \quad \text{for all } a \in A, b\in \frakB, b^{\dag} \in \frakBo.
   \end{align*}
   By \cite[Lemma 3.9]{timmer:cpmu}, the pair
   \begin{align*}
(\hba,A) \fibreh (\kdg,B):= \big( \hba \rtensorh \kdg,\AfibreB\big)    
   \end{align*}
   is a concrete $C^{*}$-$\cbasesb$-$\cbasesb$-algebra .

   For $i=1,2$, let $\phi^{(i)}$ be a morphism of nondegenerate
   concrete $C^{*}$-$\cbasesb$-$\cbasesb$-algebras
   $(\hbak{i},A^{(i)})$ and
   $({_{\delta_{i}}K^{i}_{\gamma_{i}}},B^{(i)})$. If the
   fiber product $A^{(1)} \hfibre{\alpha_{1}}{\beta_{2}} A^{(2)}$ is
   nondegenerate,  then the $*$-homomorphism
   \begin{align*}
     \phi^{(1)} \fibreh \phi^{(2)} \colon A^{(1)}
     \hfibre{\alpha_{1}}{\beta_{2}} A^{(2)} \to B^{(1)}
     \hfibre{\gamma_{1}}{\delta_{2}} B^{(2)}
   \end{align*}
   is a morphism from $(\hbak{1},A^{(1)}) \fibreh (\hbak{2},A^{(2)})$
   to $({_{\delta_{1}}K^{1}_{\gamma_{1}}},B^{(1)}) \fibreh
   ({_{\delta_{2}}K^{2}_{\gamma_{2}}},B^{(2)})$ \cite[Theorem
   3.15]{timmer:cpmu}.

   We denote by $\cbalg$ the category of all nondegenerate concrete
   $C^{*}$-$\cbasesb$-$\cbasesb$-algebras, whose morphisms between
   $C^{*}$-$\cbasesb$-$\cbasesb$-algebras $(\hba,A)$ and $(\kdg,B)$
   are all nondegenerate morphisms from $(\hba,A)$ to $(\kdg,M(B))$.
   Unfortunately, the fiber product defined above does {\em not}
   induce a natural monoidal structure on this category. Apart from
   the problem that the fiber product of nondegenerate
   $C^{*}$-$\cbasesb$-$\cbasesb$-algebras need not be nondegenerate,
   we encounter the following additional problems.

\paragraph{Unitality} 
The fiber product seems to admit a unit only on certain subcategories
of $\cbalg$.  Since the fiber product of concrete
$C^{*}$-$\cbasesb$-$\cbasesb$-algebras acts on the relative tensor
product of the underlying $C^{*}$-$\cbasesb$-$\cbasesb$-bimodules and
$\cbaseosb$ is the unit for this relative tensor product, a unit for
the fiber product should be of the form $(\cbaseosb,C)$, where $C
\subseteq {\cal L}(\frakH)$ is a suitable $C^{*}$-algebra that has to
be determined.

Given concrete $C^{*}$-$\cbasesb$-$\cbasesb$-algebras $(\cbaseosb,C)$
and $(\hba,A)$, consider the $*$-homo\-morphisms
\begin{align*}
  l_{(\hba,A)} \colon C \hfibre{\frakB}{\beta} A \to {\cal L}(H)  
\quad
\text{and}
\quad
r_{(\hba,A)} \colon A \hfibre{\alpha}{\frakBo} C \to {\cal L}(H)
\end{align*}
given by conjugation with the isomorphisms $l_{(\hba)} \colon \frakH
\htensor{\frakB}{\beta} H \to H$ and $ r_{(\hba)} \colon H
\htensor{\alpha}{\frakBo} \frakH \to H$ of Theorem
\ref{theorem:cbbimod}.
\begin{proposition} \label{proposition:cbalg-unital}
  We have
  \begin{align*}
    l_{(\hba,A)}(C \hfibre{\frakB}{\beta} A) &= \big\{ x \in {\cal
      L}(H) \,\big|\, x \rho_{\beta}(\frakB), x^{*}
    \rho_{\beta}(\frakB) \subseteq A \text{ and } x\beta, x^{*}\beta
    \subseteq [\beta C] \big\}, \\
    r_{(\hba,A)}(A \hfibre{\alpha}{\frakBo} C) &= \big\{ x \in {\cal
      L}(H) \,\big|\, x \rho_{\alpha}(\frakBo), x^{*}
    \rho_{\alpha}(\frakBo) \subseteq A \text{ and } x\alpha, x^{*}\alpha
    \subseteq [\alpha C] \big\}.
  \end{align*}
\end{proposition}
\begin{proof}
  This follows directly from the definitions and the fact that
  \begin{align*}
    l_{(\hba)}\kbeta{2} &= \beta, & l_{(\hba)} |\frakB\rangle_{1}
    &=\rho_{\beta}(\frakB), & r_{(\hba)} \kalpha{2} &= \alpha, &
    r_{(\hba)} |\frakBo\rangle_{2} &= \rho_{\alpha}(\frakBo),
  \end{align*}
  see also Remark \ref{remark:cbbimod}.
\end{proof}
\begin{remarks}
  Let $(\hba,A)$ and $(\kdg,B)$ be concrete
  $C^{*}$-$\cbasesb$-$\cbasesb$-algebras.
  \begin{enumerate}
  \item Assume that $\frakB$ and $\frakBo$ are unital. Then by
    Proposition \ref{proposition:cbalg-unital}, 
    \begin{gather*}
      l_{(\hba,A)}(\frakBo \hfibre{\frakB}{\beta} A) = A \cap {\cal
        L}(H_{\beta}) \quad \text{and} \quad r_{(\hba,A)}(A
      \hfibre{\alpha}{\frakBo} \frakB) = A \cap {\cal L}(H_{\alpha}).
    \end{gather*}
  \item Proposition \ref{proposition:cbalg-unital} suggests to
    consider the spaces
    \begin{align*}
      \tilde A:= \{ x \in {\cal L}(H) \mid x\rho_{\alpha}(\frakBo),
      x^{*}\rho_{\alpha}(\frakBo), x\rho_{\beta}(\frakB),
      x^{*}\rho_{\beta}(\frakB) \subseteq A\}, \\
      \tilde B := \{ y \in{\cal L}(K) \mid y\rho_{\gamma}(\frakBo),
      y^{*}\rho_{\gamma}(\frakBo), y\rho_{\delta}(\frakB),
      y^{*}\rho_{\delta}(\frakB) \subseteq B\}.
    \end{align*}
    Clearly, $(\hba,\tilde A)$ and $(\kdg,\tilde B)$ are
    $C^{*}$-$\cbasesb$-$\cbasesb$-algebras and $A \subseteq \tilde A$
    and $B \subseteq \tilde B$. Moreover, 
    \begin{align*}
      \tilde A \fibread \tilde B = A \fibread B
    \end{align*}
    because $[\kdelta{2}\tilde A] = [|\delta \frakBo\rangle_{2} \tilde A] =
    [\kdelta{2} \rho_{\alpha}(\frakBo)\tilde A] =  [\kdelta{2}A]$ and
    $[\kalpha{1} \tilde B] = [|\alpha \frakB\rangle_{1} \tilde B] =
    [\kalpha{1}\rho_{\delta}(\frakB) \tilde B] =
    [\kalpha{1}B]$. However, in general, $(\tilde A \fibread \tilde
    B)$ does not coincide with
    \begin{align*}
      \big\{ z \in {\cal L}(\HfibreK) \,\big|\, z\rho_{(\alpha \rt
        \gamma)}(\frakBo), z^{*}\rho_{(\alpha \rt \gamma)}(\frakBo),
      z\rho_{(\beta \lt \delta)}(\frakB), z^{*}\rho_{(\beta \lt
        \delta)}(\frakB) \subseteq \AfibreB\big\}.
    \end{align*}
  \item If $C_{1},C_{2} \subseteq {\cal L}(\frakH)$ are
    $C^{*}$-algebras and
    \begin{align*}
      A\beta &\subseteq [\beta C_{1}], &A\alpha &\subseteq [\alpha
      C_{2}], & B\delta &\subseteq [\delta C_{1}], & B\gamma
      &\subseteq [\gamma C_{2}],
    \end{align*}
  then
  \begin{align*}
    (\AfibreB) (\beta \lt \delta) &\subseteq
    [(\AfibreB)|\delta\rangle_{2}\beta] \subseteq
    [|\delta\rangle_{2}A\beta] \subseteq [|\delta\rangle_{2}\beta
    C_{1}] = [(\beta \lt \delta)C_{1}]
  \end{align*}
  and similarly $(\AfibreB) (\alpha \rt \gamma) \subseteq [(\alpha \rt
  \gamma)C_{2}]$.
  \end{enumerate}
\end{remarks}

\paragraph{Associativity} The fiber product is not associative: If
 $(\hba,A)$, $(\kdg,B)$, $(\lfe,C)$ are concrete
$C^{*}$-$\cbasesb$-$\cbasesb$-algebras, then the isomorphism
\begin{align*}
  {\cal L}\big((\HfibreK) \htensor{(\alpha \rt \gamma)}{\phi} L \big)
  \to
  {\cal L}\big(H \htensor{\alpha}{(\delta \lt \phi)} (K
  \htensor{\gamma}{\phi} L)\big)
\end{align*}
given by conjugation with the isomorphism
\begin{align*}
  \alpha_{(\hba,\kdg,\lfe)} \colon (\HfibreK) \htensor{(\alpha \rt
    \gamma)}{\phi} L \to H \htensor{\alpha}{(\delta \lt \phi)}(K
  \htensor{\gamma}{\phi} L)
\end{align*}
of Theorem \ref{theorem:cbbimod} need {\em not} identify $(\AfibreB)
\hfibre{(\alpha \rt \gamma)}{\phi} C$ with $A \hfibre{\alpha}{(\delta
  \lt \phi)} (B \hfibre{\gamma}{\phi} C)$. 

However, for each $n \geq 2$, we can define an {\em unconditional fiber
product} of $n$ $C^{*}$-$\cbasesb$-$\cbasesb$-algebras
$(\hbak{1},A^{(1)}), \ldots, (\hbak{n},A^{(n)})$ as follows.
Since the relative tensor product of
$C^{*}$-$\cbasesb$-$\cbasesb$-bimodules is associative, we can define a
 $C^{*}$-$\cbasesb$-$\cbasesb$-bimodule
\begin{align*}
  \hba := \hbak{1} \rtensorh \cdots \rtensorh \hbak{n},
\end{align*}
neglecting the order in which the relative tensor products are formed.
For each $k \in \{1,\ldots,n-1\}$, put $\alpha^{(k)} := \alpha_{1} \rt
\cdots \rt \alpha_{k}$ and $\beta^{(k+1)} := \beta_{k+1} \lt \cdots
\lt \beta_{n}$.  Let $k \in \{1,\ldots,n\}$ and put
$\sigma_{k}:=\rho_{\beta_{k}}$ and $\rho_{k}:=\rho_{\alpha_{k}}$. If
$1 < k < n$, we can identify $H$ with
\begin{align*}
  \alpha^{(k-1)} \tr_{\sigma_{k}} H^{k} {_{\rho_{k}}\tl} \beta^{(k+1)},
\end{align*}
and define $\gamma^{(k)} \subseteq {\cal L}(H^{k},H)$ to be the closed
linear span of all operators of the form $\zeta \mapsto \xi \tr \zeta
\tl \eta$, where $\xi \in \alpha^{(k-1)}$ and $\eta \in
\beta^{(k+1)}$. We put $\gamma^{(1)}:=\beta^{(2)}$ and
$\gamma^{(n)}:=\alpha^{(n-1)}$.

 The {\em unconditional fiber product} of  $(\hbak{1},A^{(1)})$,
$\ldots$, $(\hbak{n},A^{(n)})$ is the $C^{*}$-algebra
\begin{align*}
A= A^{(1)} \hfibre{\alpha_{1}}{\beta_{2}}  \cdots
 \hfibre{\alpha_{n-1}}{\beta_{n}} A^{(n)}    := \big\{ T \in {\cal
   L}(H) \,\big|\,   T\gamma^{(k)}, T^{*}\gamma^{(k)} \subseteq
 [\gamma^{(k)}A^{(k)}] \text{ for } k=1,\ldots,n\big\}.
\end{align*}
One easily checks that $({_{\beta}H_{\alpha}},A)$ is a concrete
$C^{*}$-$\cbasesb$-$\cbasesb$-algebra.

There are many different ways to form a $C^{*}$-algebra on $H$ by
successive applications of the fiber product construction to the
$C^{*}$-$\cbasesb$-$\cbasesb$-algebras $(\hbak{1},A^{(1)})$,
$\ldots$, $(\hbak{n},A^{(n)})$.  These ways correspond bijectively
with the set ${\cal T}_{n}$ of all binary trees with precisely $n$
leaves, where the leaves represent $(\hbak{1},A^{(1)})$, $\ldots$,
$(\hbak{n},A^{(n)})$, and each internal node of the tree represents
the fiber product of the $C^{*}$-$\cbasesb$-$\cbasesb$-algebras
associated to the left and to the right subtree of that node. For each
tree $t \in {\cal T}_{n}$, we denote the corresponding iterated fiber
product of $(\hbak{1},A^{(1)})$,
$\ldots$, $(\hbak{n},A^{(n)})$ by
\begin{align*}
  \bigstar^{t} \big( (\hbak{1},A^{(1)}), \ldots,
  (\hbak{n},A^{(n)}) \big) \subseteq {\cal L}(H).
\end{align*}

Now, $A^{(1)} \hfibre{\alpha_{1}}{\beta_{2}}  \cdots
 \hfibre{\alpha_{n-1}}{\beta_{n}} A^{(n)}$ is a maximal fiber product
 in the following sense:
 \begin{proposition}
   $  \bigstar^{t} \big( (\hbak{1},A^{(1)})$, $\ldots$,
   $(\hbak{n},A^{(n)}) \big) \subseteq A^{(1)}
   \hfibre{\alpha_{1}}{\beta_{2}} \cdots
   \hfibre{\alpha_{n-1}}{\beta_{n}} A^{(n)}$   for each $t \in {\cal T}_{n}$.
 \end{proposition}
 \begin{proof}
   Straightforward.
 \end{proof}
The unconditional fiber product is functorial in the following sense:
\begin{proposition}
  For each $k=1,\ldots,n$, let $\phi^{(k)}$ be a morphism of
  $C^{*}$-$\cbasesb$-$\cbasesb$-algebras $(\hbak{k},A^{(k)})$ and
  $(\kdgk{k},B^{(k)})$.  Then there exists a unique morphism
  \begin{align*}
    \phi^{(1)} \fibreh \cdots \fibreh \phi^{(n)} \colon
    A^{(1)}\hfibre{\alpha_{1}}{\beta_{2}} \cdots
    \hfibre{\alpha_{n-1}}{\beta_{n}} A^{(n)} \to B^{(1)}
    \hfibre{\gamma_{1}}{\delta_{2}} \cdots
    \hfibre{\gamma_{n-1}}{\delta_{n}} B^{(n)}
  \end{align*}
  such that for all $T_{1} \in I_{\phi^{(1)}}$, $\ldots$, $T_{n} \in
  I_{\phi^{(n)}}$ and $x \in A^{(1)}\hfibre{\alpha_{1}}{\beta_{2}}
  \cdots \hfibre{\alpha_{n-1}}{\beta_{n}} A^{(n)}$,
  \begin{align*}
    (\phi^{(1)} \fibreh \cdots \fibreh \phi^{(n)})(x) \cdot (T_{1}
    \rtensorh \cdots \rtensorh T_{n}) = (T_{1} \rtensorh \cdots
    \rtensorh T_{n}) \cdot x.
  \end{align*}
\end{proposition}
\begin{proof}
  The proof is essentially the same as in the case where $n=2$, see
  \cite[Proposition 3.13]{timmer:cpmu}.
\end{proof}

\subsection{A functor from $C^{*}$-families to concrete
  $C^{*}$-algebras}

\label{subsection:alg-functor}

As in Subsection \ref{subsection:bimod-functor},
we fix the following data:
\begin{itemize}
\item a decomposable $C^{*}$-algebra $B$ with an admissible inverse
  semigroup $\Theta \subseteq \paut(B)$,
\item a covariant representation $(\pi,\upsilon)$ of $(B,\Theta)$ on a
  Hilbert space $\frakK$ (see Section 3), where $\pi$ is faithful.
\end{itemize}
We define $\funct{}$ and $\functs{}$ as in Section 3, and put
  \begin{align*}
    \frakH &:= \funct{B}, & \frakB&:=\functs{\lft{\hmg{B}}}
    \subseteq {\cal L}(\frakH), & \frakBo &:= \functs{\rgt{\hmg{B}}}
    \subseteq {\cal L}(\frakH).
  \end{align*}
Then $\cbasesb$ is a $C^{*}$-base (Proposition  \ref{proposition:cbase}).
  We construct a  faithful functor from the category of
  $\Theta$-admissible $C^{*}$-families on $C^{*}$-$B$-bimodules to the
  category of nondegenerate concrete
  $C^{*}$-$\cbasesb$-$\cbasesb$-algebras.

\begin{proposition}
  Let $\msc$ be a $C^{*}$-family on a right $C^{*}$-$B$-$B$-bimodule
  $E$.
  \begin{enumerate}
  \item  $\functs{\msc} \subseteq {\cal L}(\funct{E})$ is
  $C^{*}$-algebra, and $\functs{\msm(\msc)} \subseteq M(\functs{\msc})$.
\item If $\msc$ is $\Theta$-admissible, then
  $\big({_{\beta(E)}\funct{(E)}_{\alpha(E)}},\, \functs{\msc}\big)$ is
  a nondegenerate concrete $C^{*}$-$\cbasesb$-$\cbasesb$-algebra.
\end{enumerate}
\end{proposition}
\begin{proof}
  i) By Proposition \ref{proposition:functs},
  $\functs{(\msc)}^{*}=\functs{(\msc^{*})}=\functs{\msc}$ and
  $[\functs{(\msc)}\functs{(\msc)}] \subseteq \functs{[\msc\msc]}=
  \functs{\msc}$, so $\functs{\msc}$ is a $C^{*}$-algebra. Likewise,
  $\functs{\msm(\msc)}$ is a $C^{*}$-algebra, and
  $\functs{\msm(\msc)} \subseteq M(\functs{\msc})$ because
  $[\functs{(\msm(\msc))} \functs{(\msc)}] \subseteq
  \functs{[\msm(\msc)\msc]} = \functs{\msc}$.

\smallskip

ii)    By Propositions \ref{proposition:functs}, \ref{proposition:modules}
  and $\Theta$-admissibility of $\msc$,
  \begin{align*}
    [\rho_{\alpha(E)}(\frakBo)\functs{\msc}] &=
    \functs{[\rgte{\hmg{B}}\msc]} = \functs{\msc}, &
    [\rho_{\beta(E)}(\frakB) \functs{\msc}] &= \functs{[\lfte{\hmg{B}}
      \msc]} = \functs{\msc}.
  \end{align*}
  By \cite[Remarks 3.9 (iv)]{timmermann:hopf}, the
  $C^{*}$-family $\msc$ is nondegenerate if and only if the
  $C^{*}$-algebra $\msc^{\Id}_{\Id} \subseteq {\cal L}^{B}_{B}(E)$ is
  nondegenerate, and in that case, also $\functs{\msc^{\Id}_{\Id}}
  \subseteq \functs{\msc} \subseteq {\cal L}(\funct{E})$ is
  nondegenerate.  
\end{proof}
The assignment $\msc \mapsto \functs{\msc}$ is
functorial in the following sense:
\begin{proposition}
  Let $\msc$ and $\msd$ be $\Theta$-admissible $C^{*}$-families on
  right $C^{*}$-$B$-$B$-bimodules  $E$ and $F$,
  respectively.
  \begin{enumerate}
  \item Let $\phi$ be a normal  morphism from $\msc$ to
    $\msm(\msd)$. Then there exists a unique 
    $*$-homo\-morphism $\functs \phi \colon \functs \msc \to M(\functs
    \msd)$ such that $(\functs{\phi})( \funct{c}) =
    \funct{\phi^{\rho}_{\sigma}(c)}$ for all $c \in
    \msc^{\rho}_{\sigma}$, $\rho,\sigma \in \Theta$.  If $\phi$ is nondegenerate, then so is
    $\functs{\phi}$. If $\phi(\msc) \subseteq \msd$, then
    $(\functs{\phi})(\functs{\msc}) \subseteq \functs{\msd}$. Moreover,
    $\functs{\phi}$ is a morphism from
    $\big({_{\beta(E)}\funct{(E)}_{\alpha(E)}},\, \functs{\msc}\big)$
    to $\big({_{\beta(F)}\funct{(F)}_{\alpha(F)}},\,
    M(\functs{\msd})\big)$. 
  \item Let $\psi$ be a morphism from
    $\big({_{\beta(E)}\funct{(E)}_{\alpha(E)}},\, \functs{\msc}\big)$
    to $\big({_{\beta(F)}\funct{(F)}_{\alpha(F)}},\,
    M(\functs{\msd})\big)$.  Then there exists a $\Theta$-admissible
    $C^{*}$-family $\msb \subseteq \cL(F)$ and a unique normal morphism
    $\phi$ from $\msc$ to $\msb$ such that $\psi = \functs{\phi}$.
  \end{enumerate}
\end{proposition}
\begin{proof}
  i) Existence of the $*$-homomorphism $\functs{\phi}$ follows by a
  similar argument as in the proof of Proposition
  \ref{proposition:morphism-extend}.  We only  prove that $\functs{\phi}$ is
  a morphism of $C^{*}$-$\cbasesb$-$\cbasesb$-algebras,  i.e., that
  $[I_{\functs{\phi}}\alpha(E)]\supseteq \alpha(F)$ and
  $[I_{\functs{\phi}}\beta(E)]\supseteq \beta(F)$.  Using Proposition
  \ref{proposition:module-morphisms}, one easily checks that
  $\functs{I_{\phi}} \subseteq I_{\functs{\phi}}$. Now, $[I_{\phi}E]=
  F$ by assumption of $I_{\phi}$, and hence
  $[I_{\functs{\phi}}\alpha(E)] \supseteq
      \functs{[I_{\phi}\lft{\hmg{E}}]} = \functs{\lft{\hmg{F}}} = \alpha(F)$.
    A similar calculation shows that
    $[I_{\functs{\phi}}\alpha(E)]\supseteq \beta(F)$. 

  \smallskip

  ii) Put $J_{\psi}:=\{ T \in {\cal L}^{B}_{B}(E,F) \mid \funct T \in
  I_{\psi}\}$. Then $[(\functs{J_{\psi}})\alpha(E)] =
  [I_{\psi}\alpha(E)]=\alpha(F)$ by Proposition
  \ref{proposition:module-morphisms} and assumption on $\psi$, and
  hence $[J_{\psi}E]=F$. 

  Note that $[J_{\psi}^{*}J_{\psi}]$ commutes with every element of
  $\msc$ because $\msc$ is $\Theta$-supported,
  $[I_{\psi}^{*}I_{\psi}]=\functs{[J_{\psi}^{*}J_{\psi}]}$ commutes
  with $\functs{\msc}$, and $\funct{}$ is faithful.  Let $\rho,\sigma
  \in \paut(B)$ and $c \in \msc^{\rho}_{\sigma}$. Since $\langle
  S'\xi'|Sc\xi\rangle = \sigma(\langle S'c^{*}\xi'|S\xi\rangle)$ for
  all $S,S' \in J_{\psi}$ and $\xi,\xi' \in E$, we can define linear
  maps
  \begin{align*}
    \phi^{\rho}_{\sigma}(c)^{*} \colon F \to F, \ S'\xi' \mapsto
    S'c^{*}\xi', \quad \text{and} \quad
    \phi^{\sigma}_{\sigma}(c) \colon F \to F, \ S\xi \mapsto Sc\xi,
  \end{align*}
  where $S,S' \in J_{\psi}$ and $\xi,\xi' \in E$. Using the relation
  $J_{\psi}^{*}J_{\psi} \in {\cal L}^{B}_{B}(E)$, one easily checks
  that $\phi^{\rho}_{\sigma}(c)^{*} \in
  \cL^{\rho^{*}}_{\sigma^{*}}(F)$ and $\phi^{\rho}_{\sigma}(c) \in
  \cL^{\rho}_{\sigma}(F)$.  Letting $\rho,\sigma,c$ vary, we obtain a
  family of maps $\phi^{\rho}_{\sigma} \colon \msc^{\rho}_{\sigma} \to
  \cL^{\rho}_{\sigma}(F)$, where $\rho,\sigma \in \Theta$. By
  construction, $J_{\psi} \subseteq I_{\phi}$. Thus, $\phi$ is a
  normal morphism from $\msc$ to $\cL(F)$.  By Proposition
  \ref{proposition:family-morphism} iv), 
  $\msb:=([\phi^{\rho}_{\sigma}(\msc^{\rho}_{\sigma})])_{\rho,\sigma}
  \subseteq \cL(F)$ is a $\Theta$-admissible $C^{*}$-family. By
  construction, $\psi=\functs{\phi}$.
\end{proof}
\begin{remark}
  We do not know whether the normal morphism $\phi$ in ii) 
  satisfies $\phi(\msc) \subseteq \msd$.
\end{remark}

\begin{corollary}
  The assignments $\msc \mapsto \functs{\msc}$ and $\phi \mapsto
  \functs{\phi}$ define a faithful functor
  \begin{align*}
     \bfamily \to \cbalg.
  \end{align*}
\end{corollary}

The functor constructed above embeds the fiber product of
$C^{*}$-families into the spatial $C^{*}$-fiber product of concrete
$C^{*}$-$\cbasesb$-$\cbasesb$-algebras. Indeed, let $\msc$ and $\msd$
be $\Theta$-admissible $C^{*}$-families on right
$C^{*}$-$B$-$B$-bimodules $E$ and $F$, respectively. Then conjugation by
the isomorphism $\tau_{E,F} \colon \funct{(E)}
\htensor{\alpha(E)}{\beta(F)} \funct{(F)} \to \funct{(E \tr F)}$
defines an isomorphism
\begin{align*}
  \Ad_{\tau_{E,F}} \colon {\cal L}\big(\funct{(E)}
  \htensor{\alpha(E)}{\beta(F)} \funct{(F)}\big) \to
  {\cal L}(\funct{(E \tr F)}),
\end{align*}
and the following result holds:
\begin{theorem} \label{theorem:fibre-embed}
  Let $\msc$ and $\msd$ be $\Theta$-admissible $C^{*}$-families on
 right $C^{*}$-$B$-$B$-bimodules $E$ and $F$,
  respectively. Then
  $\Ad_{\tau_{E,F}}^{-1}\big( \functs{(\msc \tr \msd)}\big) \subseteq 
    \functs{(\msc)} \hfibre{\alpha(E)}{\beta(F)} \functs{(\msd)}$.
  \end{theorem}
\begin{proof}
  Put $\tau:=\tau_{E,F}$.  It suffices to
  show that
  \begin{align*}
    \big[\functs{(\msc \tr \msd)}\tau|\alpha(E)\rangle_{1}\big] \subseteq
    \big[\tau|\alpha(E)\rangle_{1} \functs{(\msd)}\big] \quad \text{and} \quad
    \big[\functs{(\msc \tr \msd)}\tau|\beta(F)\rangle_{2}\big] \subseteq
    \big[\tau|\beta(F)\rangle_{2} \functs{(\msc)}\big].
  \end{align*}
  Since $\tau|\alpha(E)\rangle_{1}= \functs{|\hmg{E}\rangle_{1}}$ and
  $\tau|\beta(F)\rangle_{2} = \functs{|\hmg{F}]_{2}}$ (see
equation  \eqref{eq:tau-ketbra}), these inclusions follow from the relations
  \begin{align*}
   \functs{[(\msc \tr \msd) |\hmg{E}\rangle_{1}]} &=
  \functs{[|\msc\hmg{E}\rangle_{1}\msd]} =
  \functs{[|\hmg{E}\rangle_{1}\msd]}, \\
  \functs{[(\msc \tr \msd)|\hmg{F}]_{2}]} &=
  \functs{[|\msd\hmg{F}\rangle_{1}\msc]} =
  \functs{[|\hmg{F}\rangle_{2}\msc]}. \qedhere
  \end{align*}
\end{proof}
\begin{theorem} \label{theorem:cfam-morphism-tensor} For $i=1,2$, let
  $\msc^{(i)}$ and $\msd^{(i)}$ be $\Theta$-admissible
  $C^{*}$-families on $C^{*}$-$B$-$B$-bimodules $E^{i}$ and
  $F^{i}$, respectively, and let $\phi^{(i)}$ be a nondegenerate
  normal morphism from $\msc^{(i)}$ to $\msm(\msd^{(i)})$.  Put
  $\Phi_{E} := \tau_{E^{1},E^{2}}$ and
  $\Phi_{F}:=\tau_{F^{1},F^{2}}$.  Then the following diagram
  commutes:
  \begin{align*}
    \xymatrix@C=60pt{ {\functs{(\msc^{(1)} \tr \msc^{(2)})}}
      \ar[r]^{\functs{(\phi^{(1)}\tr \phi^{(2)})}
      } \ar@{^(->}[d]^{\Ad^{-1}_{\Phi_{E}}}& {
        M\big(\functs{(\msd^{(1)}
          \tr \msd^{(2)})}\big)}
      \ar@{^(->}[d]^{\Ad^{-1}_{\Phi_{F}}}\\
      {(\functs{\msc^{(1)}})
        \hfibre{\alpha(E^{1})}{\beta(E^{2})}
        (\functs{\msc^{(2)}}})  \ar[r]^{\functs{(\phi^{(1)})}
        \fibreh \functs{(\phi^{(2)})}} &
      {M\big((\functs{\msd^{(1)}})
      \hfibre{\alpha(F^{1})}{\beta(F^{2})}
      (\functs{\msd^{(2)}})\big). }}
  \end{align*}
\end{theorem}
\begin{proof}
  This follows easily from the definitions.
\end{proof}

\section{Applications to  pseudo-multiplicative unitaries and Hopf $C^{*}$-bimodules}

In this section, we apply the functors constructed in Subsections
\ref{subsection:bimod-functor} and \ref{subsection:alg-functor} to
Hopf $C^{*}$-families and to pseudo-multiplicative unitaries.  As
before, we fix a decomposable $C^{*}$-algebra $B$ with an admissible
inverse semigroup $\Theta \subseteq \paut(B)$ and a covariant
representation $(\pi,\upsilon)$ of $(B,\Theta)$ on a Hilbert space
$\frakK$ (see Section 3), where $\pi$ is faithful, define $\funct{}$,
$\functs{}$ as in Section 3, and put $\frakH := \funct{B}$, $\frakB
:=\functs{\lft{\hmg{B}}} \subseteq {\cal L}(\frakH)$, $\frakBo :=
\functs{\rgt{\hmg{B}}} \subseteq {\cal L}(\frakH)$.

\subsection{From Hopf $C^{*}$-families to concrete Hopf
  $C^{*}$-bimodules}

Recall that
a {\em Hopf $C ^{*}$-family} \cite{timmermann:hopf} over $B$ is a
nondegenerate $C^{*}$-family $\msa$ on a right
$C^{*}$-$B$-$B$-bimodule equipped with a nondegenerate morphism
$\Delta \colon \msa \to \msm(\msa \tr \msa)$ such that
  \begin{enumerate}
  \item $[\Delta(\msa)(\Id \tr \msa^{\Id})] \subseteq \msa \tr
    \msa$ and $[\Delta(\msa)(\msa_{\Id} \tr \Id)] \subseteq
    \msa \tr \msa$, and
  \item the following diagram commutes:
    \begin{align*}
      \xymatrix@R=15pt@C=35pt{ {\msa} \ar[r]^(0.45){\Delta}
        \ar[d]^{\Delta} & {\msm(\msa \tr \msa)}
        \ar[d]^{\Id \tr \Delta} \\
        {\msm(\msa \tr \msa)} \ar[r]^(0.45){\Delta \tr
          \Id} & {\msm(\msa \tr \msa \tr \msa).}  }
    \end{align*}
  \end{enumerate}
  We call $(\msa,\Delta)$ {\em bisimplifiable} if the inclusions in i)
  are equalities, and {\em $\Theta$-admissible} if $\msa$ is
  $\Theta$-admissible and $\Delta$ is normal.  A normal morphism of
  $\Theta$-admissible Hopf $C^{*}$-families $(\msa,\Delta_{\msa})$ and
  $(\msc,\Delta_{\msc})$ is a normal nondegenerate morphism $\phi$
  from $\msa$ to $\msm(\msc)$ that makes the following diagram
  commute:
    \begin{align*}
      \xymatrix@R=15pt@C=35pt{
      \msa \ar[r]^{\phi} \ar[d]^{\Delta_{\msa}} & \msm(\msc)
      \ar[d]^{\Delta_{\msc}} \\
      \msm(\msa \tr \msa) \ar[r]^{\phi \tr \phi} & 
      \msm(\msc \tr \msc).
      }
    \end{align*}
    Replacing the internal tensor product ``$\tr$'' by the flipped
    internal tensor product ``$\tl$'' in the target of $\Delta$ and in
    the diagrams above, one arrives at the notion of  {\em flipped
      Hopf $C^{*}$-families} and their morphisms. 

 Recall that a {\em concrete Hopf $C^{*}$-bimodule}
    \cite{timmer:cpmu} over $\cbasesb$ consists of a nondegenerate
    concrete $C^{*}$-$\cbasesb$-$\cbasesb$-algebra $(\hba,A)$ and a
    morphism $\Delta$ from $(\hba,A)$ to $(\hba,A) \fibreh (\hba,A)$
    that makes the following diagram commute,
    \begin{align*}
      \xymatrix@C=35pt@R=15pt{ {A} \ar[r]^(0.45){\Delta}
        \ar[d]^{\Delta} & {A \hfibre{\alpha}{\beta} A}
        \ar[d]^{\Id \fibreh \Delta} \\
        {A \hfibre{\alpha}{\beta} A} \ar[r]^(0.45){\Delta
          \fibreh \Id} & {A \hfibre{\alpha}{\beta} A
          \hfibre{\alpha}{\beta} A,}  }
    \end{align*}
    where $A \hfibre{\alpha}{\beta} A \hfibre{\alpha}{\beta} A$ denotes
    the unconditional fiber product (see Subsection \ref{subsection:cbalg}).
  A morphism of concrete Hopf $C^{*}$-bimodules $(\hba,A,\Delta_{A})$
  and $(\kdg,C,\Delta_{C})$ is a  morphism $\phi$ from
  $(\hba,A)$ to $(\kdg,M(C))$ that makes the following diagram
  commute:
   \begin{align*}
     \xymatrix@R=15pt@C=35pt{
       A \ar[r]^{\phi} \ar[d]^{\Delta_{A}} & M(C) \ar[d]^{\Delta_{C}} \\
       A \hfibre{\alpha}{\beta} A \ar[r]^{\phi \fibreh \phi} &
       M(C \hfibre{\gamma}{\delta} C).
     }
   \end{align*}
  \begin{remark}
    We can not yet formulate an analogue of the bisimplifiability
    condition for a concrete Hopf $C^{*}$-bimodule $(\hba,A,\Delta)$.
    One could compare spaces of the form $[\Delta(A)(A^{(\alpha)}
    \htensor{\alpha}{\beta} 1)]$ and $[\Delta(A)(1
    \htensor{\alpha}{\beta} A^{(\beta)})]$, where natural choices for
    $A^{(\alpha)}$ and $A^{(\beta)}$ are $A \cap
    \rho_{\alpha}(\frakBo)'$ and $A \cap \rho_{\beta}(\frakB)'$,
    respectively, or $A \cap {\cal L}(H_{\alpha})$ and $A \cap {\cal
      L}(H_{\beta})$, to the fiber product $A \fibreab A$. But this
    fiber product has to be replaced by a smaller $C^{*}$-subalgebra.

    Another natural condition on a concrete Hopf $C^{*}$-bimodule
    $(\hba,A,\Delta)$ would be to demand that
    $[\Delta(A)\kalpha{1}]=[\kalpha{1}A]$ and
    $[\Delta(A)\kbeta{2}]=[\kbeta{2}A]$.
  \end{remark}

  The functor constructed in the preceding section yields an embedding
  of the category of $\Theta$-admissible Hopf $C^{*}$-families into
  the category of concrete Hopf $C^{*}$-bimodules:
  \begin{theorem} \label{theorem:hopf}
    \begin{enumerate}
    \item Let $(\msa,\Delta_{\msa})$ be a $\Theta$-admissible Hopf
      $C^{*}$-family on a right $C^{*}$-$B$-$B$-bimodule $E$.  Then
      $({_{\beta(E)}\funct{(E)}_{\alpha(E)}},\functs{\msa})$ together
      with the composition
      \begin{align*}
        \Delta_{\functs{\msa}} := \Ad_{\tau_{E,E}}^{-1} \circ
        \functs{(\Delta_{\msa})} \colon \functs{\msa} \to
        M(\functs{(\msa \tr \msa)}) \to M(\functs{(\msa)}
        \hfibre{\alpha}{\beta} \functs{(\msa)})
      \end{align*}
      is a concrete Hopf $C^{*}$-bimodule.
    \item Let $(\msa,\Delta_{\msa})$ and $(\msc,\Delta_{\msc})$ be
      $\Theta$-admissible Hopf $C^{*}$-families on right
      $C^{*}$-$B$-$B$-bimodules $E$ and $F$, respectively, and let
      $\phi$ be normal morphism from $(\msa,\Delta_{\msa})$ to
      $(\msc,\Delta_{\msc})$.  Put
      $\Delta_{\functs{\msc}}:=\Ad_{\tau_{F,F}}^{-1} \circ
      \functs{(\Delta_{\msc})}$. Then $\functs{\phi} \colon
      \functs{\msa}\to M(\functs{\msc})$ is a morphism of the concrete
      Hopf $C^{*}$-bimodules $({_{\beta(E)}\funct{(E)}_{\alpha(E)}},
      \functs{\msa}, \Delta_{\functs{\msa}})$ and
      $({_{\beta(F)}\funct{(F)}_{\alpha(F)}}, \functs{\msc},
      \Delta_{\functs{\msc}})$.
    \end{enumerate}
  \end{theorem}
  \begin{proof}
    i) Put $\alpha:=\alpha(E)$, $\beta:=\beta(E)$, and
    $A:=\functs{\msa}$. First, we show that $\Delta_{A}(A) \subseteq A
    \hfibre{\alpha}{\beta} A$.  Since $\msa$ is nondegenerate, so is
    the $C^{*}$-algebra $\msa^{\Id}_{\Id}$.  Using the relation
    $\tau_{E,E}|\alpha(E)\rangle_{1} = \functs{|\hmg{E}\rangle_{1}}$
    (see equation \eqref{eq:tau-ketbra}) and Theorem \ref{theorem:fibre-embed}
    \begin{align*}
      \big[\Delta_{A}(A) \kalpha{1}\big] &\subseteq \tau_{E,E}
      \functs{\big[\Delta_{\msa}(\msa)|\hmg{E}\rangle_{1}\big]} \\ &=
      \tau_{E,E} \functs{\big[\Delta_{\msa}(\msa)(\msa^{\Id}_{\Id} \tr
        \Id)|\hmg{E}\rangle_{1}\big]} \\
      &\subseteq \tau_{E,E} \functs{\big[(\msa \tr
        \msa)|\hmg{E}\rangle_{1}\big]} \\ &=
      \big[\Ad_{\tau_{E,E}}(\functs{(\msa \tr
        \msa)})|\alpha\rangle_{1}\big] \subseteq
      [(\AfibreA)\kalpha{1}] \subseteq [|\alpha\rangle_{1} A],
    \end{align*}
    and a similar calculation shows that $[\Delta_{A}(A)\kbeta{2}]
    \subseteq [\kbeta{2}A]$. Therefore, $\Delta_{A}(A) \subseteq A
    \hfibre{\alpha}{\beta} A$.

    The equation $(\Delta_{A} \fibreh \Id) \circ \Delta_{A} = (\Id
    \fibreh \Delta_{A}) \circ \Delta_{A}$ follows from the equation
    $(\Delta_{\msa} \tr \Id) \circ \Delta_{\msa} = (\Id
    \tr \Delta_{\msa}) \circ \Delta_{\msa}$ and Theorem
    \ref{theorem:cfam-morphism-tensor}.

    \smallskip
    
    ii) This follows directly from Theorem
    \ref{theorem:cfam-morphism-tensor}. 
  \end{proof}

  \begin{remark} \label{remark:hopf-flipped}
    The constructions in the preceding theorem carry over to flipped
    Hopf $C^{*}$-families in a straightforward way.
  \end{remark}

\subsection{Pseudo-multiplicative unitaries on $C^{*}$-modules}

Finally, we apply the techniques developed so far to
pseudo-multiplicative unitaries, and compare the approaches of
\cite{timmermann,timmermann:hopf} and \cite{timmer:ckac,timmer:cpmu}.

Let us recall the definition of  pseudo-multiplicative unitaries on
$C^{*}$-modules from \cite{timmermann,timmermann:hopf}.  A {\em
  $C^{*}$-trimodule} over $B$ consists of a full $C^{*}$-module $E$
over $B$ and two faithful nondegenerate commuting representations
$s,r\colon B \to {\cal L}_{B}(E)$.  We denote by $_{r}E$ and $_{s}E$
the right $C^{*}$-$B$-$B$-bimodules formed by $E$ and the
representations $r$ and $s$, respectively.  We call $(E,s,r)$
{$\Theta$-admissible} if $_{r}E$ and $_{s}E$ are $\Theta$-admissible.
Let $(E,s,r)$ be a $C^{*}$-trimodule over $B$. Then we can define
representations $r_{1},s_{2}, r_{2}$ on $E {_{s}\tl} E$ by
$r_{1}(b):=r(b) \tl 1$, $s_{2}(b):= 1 \tl s(b)$,
$r_{2}(b):=1 \tl r(b)$ for all $b \in B$, and similarly
representations $r_{1},s_{1}, s_{2}$ on $E \tr_{r} E$. A
{\em pseudo-multiplicative unitary on $(E,s,r)$} is a unitary 
\begin{align*}
 W
\colon E {_{s} \tl} E \to E \tr_{r} E 
\end{align*}
that satisfies the
following conditions:
\begin{enumerate}
\item $Wr_{2}(b) = s_{1}(b)W$, $Wr_{1}(b) = r_{1}(b)W$, $Ws_{2}(b) =
  s_{2}(b)W$ for all $b \in B$,
\item  the following diagram commutes,
  \begin{gather} \label{eq:pentagon-w}
    \xymatrix{ {E {_{s} \tl} E {_{s}\tl} E} \ar[r]^{W \tl 1}
      \ar[d]^{1 \tl W} & {E \tr_{r} E {_{s}\tl} E} \ar[r]^{1 \tr
        W} &
      {E \tr_{r} E \tr_{r} E,} \\
      {E {_{s}\tl} (E \tr_{r} E)} \ar[rr]^{W_{13}} && {(E {_{\hat
            r}\tl} E) \tr_{r} E} \ar[u]^{W \tr 1} }
  \end{gather}
  where $W_{13}$ acts like $W$ on the first and third component of the
  internal tensor product.
\end{enumerate}
If $W$ is such a pseudo-multiplicative unitary, then the unitary
$W^{op}:= \Sigma W^{*} \Sigma \colon E {_{r}\tl} E \to E \tr_{s}
E$ is a pseudo-multiplicative unitary on $(E,r,s)$ \cite[Remark
2.4]{timmermann:hopf}.

Next, we recall the definition of $C^{*}$-pseudo-multiplicative
unitaries from \cite{timmer:ckac,timmer:cpmu}. A {\em
  $C^{*}$-$\cbasesb$-$\cbaseosb$-$\cbaseosb$-module} consists of a
Hilbert space $H$ with pairwise compatible $C^{*}$-factorizations
$\alpha \in \cfact(H;\cbasesb)$ and $\beta,\hat\beta \in
\cfact(H;\cbaseosb)$. 
A {\em $C^{*}$-pseudo-multiplicative unitary} on a
$C^{*}$-$\cbasesb$-$\cbaseosb$-$\cbaseosb$-module
$(H,\alpha,\hat\beta,\beta)$ is a unitary
\begin{align*}
  V \colon \Hsource \to \Hrange
\end{align*}
that satisfies the following conditions:
\begin{enumerate}
\item $V(\alpha \lt \alpha) = \alpha \rt \alpha$, $V(\hat\beta \rt
  \beta) = \hat\beta \lt \beta$, $V(\hat \beta \tr\hat\beta) = \alpha
  \rt \hat \beta$, $V(\beta \lt \alpha) = \beta \lt \beta$, 
\item the following diagram commutes,
  \begin{align} \label{eq:pentagon-v}
    \xymatrix{
      {\Hone} \ar[d]^{1 \rtensorh V} \ar[r]^{V \rtensorh 1} & {\Htwo}
      \ar[r]^{1 \rtensorh V} & {\Hthree,}\\
      {\Hfive} \ar[rr]^{V_{13}} && {\Hfour} \ar[u]^{V \rtensorh 1}
    }
  \end{align}
  where $V_{13}$ acts like $W$ on the first and third component of the
  internal tensor product (see \cite[Lemma 4.1]{timmer:cpmu}).
\end{enumerate} 
If $V$ is such a $C^{*}$-pseudo-multiplicative unitary, then  the
unitary $V^{op}:=\Sigma V^{*} \Sigma \colon H \htensor{\beta}{\alpha}
H \to H \htensor{\alpha}{\hat\beta} H$ is a
$C^{*}$-pseudo-multiplicative unitary on $(H,\alpha,\beta,\hat\beta)$.
\begin{theorem}
  Let $(E,s,r)$ be a $\Theta$-admissible $C^{*}$-trimodule over
  $B$ and let $W \colon E {_{s}\tl} E \to E \tr_{r} E$ be a
  pseudo-multiplicative unitary.  Then $\alpha(\sE)=\alpha(\rE)$. Put
  \begin{gather*}
    \begin{aligned}
      H&:= \funct{E}, & \alpha &:=
      \alpha(\rE)=\alpha(\sE), & \hat\beta &:= \beta(\sE), & \beta &:=
      \beta(\rE),
    \end{aligned}\\
    \begin{aligned}
      \tau &:= \tau_{\rE,\rE} \colon H \htensor{\alpha}{\beta} H \to
      \funct{(E \tr_{r} E)}, & \tau^{op} &:= \tau_{\sE,\sE} \colon H
      \htensor{\hbeta}{\alpha} H \to \funct{(E {_{s} \tl} E)}.
    \end{aligned}
  \end{gather*}
  Then $(H,\alpha,\hbeta,\beta)$ is a
  $C^{*}$-$\cbasesb$-$\cbaseosb$-$\cbaseosb$-module, and the unitary
  \begin{align*}
    V:= \tau^{-1} \circ \funct{W} \circ \tau^{op} \colon
    \Hsource \to \funct{(E {_{s}\tl} E)} \to \funct{(E \tr_{r}
    E)} \to \Hrange
  \end{align*}
  is a $C^{*}$-pseudo-multiplicative unitary.
\end{theorem}
\begin{proof}
  Remark \ref{remark:alpha} shows that $\alpha(\sE)=\alpha(\rE)$.  By
  equation \eqref{eq:tau-ketbra}, Proposition
\ref{proposition:itp},  and assumption on $W$, 
  \begin{align*}
    \tau V [\khbeta{1}\beta] =
    \funct{(W)}\functs{\big[|\sE]_{1}\rgt{\rE}\big]} &=
    \functs{\big[W\rgt{\sE \tl \rE}]} \\ &= \functs{\rgt{\sE \tr \rE}} =
    \functs{\big[|\rE]_{2}\rgt{\sE}\big]]}= \tau [\kbeta{1}\hbeta],
  \end{align*}
  showing that $V(\hbeta \tr \beta)=\hbeta \tl \beta$. Similar
  calculations show that $V(\alpha \lt \alpha) = \alpha \rt \alpha$,
  $V(\hat \beta \tr\hat\beta) = \alpha \rt \hat \beta$, $V(\beta \lt
  \alpha) = \beta \lt \beta$; here, one has to use Proposition
  \ref{proposition:ket-bra}.  Finally, tedious but straightforward
  arguments show that commutativity of diagram \eqref{eq:pentagon-w}
  implies commutativity of  diagram \eqref{eq:pentagon-v}.
\end{proof}

Let $E,s,r,W$ and $H,\alpha,\hat\beta,\beta, V$ be as in
the Theorem above. In \cite{timmermann:hopf}, we associated to $W$ two
families 
\begin{align*}
  \hmsa &:= [[\hmg{_{r}E}|_{2}W|\hmg{_{r}E}\rangle_{2}] \subseteq
  \cL(\sE) & \msa &:= [\langle \hmg{_{s}E}|_{1}W|\hmg{_{s}E}]_{1}]
  \subseteq \cL(\rE)
\end{align*}
and two families of maps 
\begin{align*} 
  (\hDelta_{W})^{\rho}_{\sigma} \colon \hmsa^{\rho}_{\sigma} &\to
  \cL^{\rho}_{\sigma}(\Esource), & (\Delta_{W})^{\rho}_{\sigma} \colon
  \msa^{\rho}_{\sigma} &\to \cL^{\rho}_{\sigma}(\Erange),
  \\
  \hat a &\mapsto W^{*}(\Id \tr \hat a)W, & a &\mapsto W(a \tl
  \Id)W^{*},
\end{align*}
where $\rho,\sigma \in \paut(B)$.  In \cite{timmer:cpmu}, we
associated to $V$ two spaces
\begin{align*}
  \hA &:= [\bbeta{2}V\kalpha{2}] \subseteq {\cal L}(H), &
  A &:= [\balpha{1}V\khbeta{1}] \subseteq {\cal L}(H),
\end{align*}
and two maps
\begin{align*}
  \hDelta_{V} \colon \hA &\to {\cal L}(\Hsource), &
  \Delta_{V} \colon  A &\to {\cal L}(\Hrange), \\
  x &\mapsto V^{*}(1 \htensor{\alpha}{\beta} x)V, & y &\mapsto V(y
  \htensor{\hbeta}{\alpha} 1)V^{*}.
\end{align*}
Recall that in Theorem \ref{theorem:hopf} and Remark
\ref{remark:hopf-flipped}, we associated to every (flipped)
$\Theta$-admissible Hopf $C^{*}$-family a concrete Hopf
$C^{*}$-bimodule.
\begin{proposition}
  \begin{enumerate}
  \item If $(\msa,\Delta_{W})$ is a Hopf $C^{*}$-family, then it is
    $\Theta$-admissible, and then $(\hba,A,\Delta_{V})$ is the associated
    concrete Hopf $C^{*}$-bimodule.
  \item If $(\hmsa,\hDelta_{W})$ is a flipped Hopf $C^{*}$-family, then it is
    $\Theta$-admissible, and $({_{\alpha}H_{\hbeta}},\hA,\hDelta_{V})$
    is the associated concrete Hopf $C^{*}$-bimodule.
  \end{enumerate}
\end{proposition}
\begin{proof}
  i) Assume that $\msa$ is a $C^{*}$-family. The definition of
  $\msa$ and the fact that $_{s}E$ is $\Theta$-decomposable imply
  that $\msa$ is $\Theta$-supported. By \cite[Proposition
  4.4]{timmermann:hopf}, $[\msa \mso(_{r}E)] \subseteq [\msa]$ and
  $[\msa^{\rho}_{\sigma}\mso^{\rho^{*}\rho}_{\sigma^{*}\sigma}({_{r}E})]
  = \msa^{\rho}_{\sigma}$ for all $\rho,\sigma \in \paut(B)$, and by
  \cite[Proposition 4.5]{timmermann:hopf}, $\msa$ is
  nondegenerate. Consequently, $\msa$ is $\Theta$-admissible.
 By definition, equation \eqref{eq:tau-ketbra}, and
  Proposition \ref{proposition:functs},
  \begin{align*} 
    A = \big[\balpha{1}\tau (\funct{W})\tau^{op} \khbeta{1}\big] =
    \functs{\big[\langle\hmg{_{s}E}|_{1}W|\hmg{_{s}E}]_{1}\big]} =
    \functs{\msa},
  \end{align*}
  and $\Delta_{V} (\funct{ a}) = \Ad_{\tau}
  \big(\funct{(\Delta_{W})^{\rho}_{\sigma}(a)}\big)$ for all
  $\rho,\sigma \in \Theta$ and $a \in \msa^{\rho}_{\sigma}$. The
  claims follow.

\smallskip

  ii) The proof is similar to the proof of i).
\end{proof}

In \cite{timmermann} and \cite{timmer:cpmu}, we studied
regularity conditions on pseudo-multiplicative unitaries that ensure
that $(\hmsa,\hDelta_{W})$ and $(\msa,\Delta_{W})$ are Hopf
$C^{*}$-families, and that $({_{\alpha}H_{\hat\beta}},\hA,\hDelta_{V})$
$(\hba,A,\Delta_{V})$ are concrete Hopf $C^{*}$-bimodules: $W$ is {\em
  regular} \cite{timmermann} if $[\langle E|_{1}W|E\rangle_{2}]
= {\cal K}_{B}(E)$, and $V$ {\em regular} \cite{timmer:cpmu} if
$[\balpha{1}V\kalpha{1}]=[\alpha\alpha^{*}]$.
\begin{proposition}
  If $W$ is regular, then $V$ is regular.
\end{proposition}
\begin{proof}
  This follows from the relations
  \begin{align*}
    [\alpha\alpha^{*}] &=
    \functs{[\lft{\hmg{_{r}E}}\lft{\hmg{_{r}E}}^{*}]} = \functs{{\cal
        K}_{B}(E)}, & [\balpha{1} V \kalpha{2}] &= \functs{[\langle
      E|_{1}W|E\rangle_{2}]}. \qedhere
  \end{align*}
\end{proof}

\def\cprime{$'$}

\end{document}